\documentclass[12pt]{amsart}
\usepackage{latexsym,fancyhdr,amssymb,color,amsmath,amsthm,graphicx,listings,comment}
\usepackage[section]{placeins}
\newtheorem{thm}{Theorem}
\newtheorem{lemma}{Lemma}

\let\paragraph\subsection
\title{Morse and Lusternik-Schnirelmann for graphs}
\author{Oliver Knill}
\date{Mai 29, 2024}
\address{Department of Mathematics \\ Harvard University \\ Cambridge, MA, 02138 }
\makeindex


\begin{document}
\maketitle

\begin{abstract}
Both Morse theory and Lusternik-Schnirelmann theory link algebra, topology and analysis
in a geometric setting. The two theories can be formulated in finite geometries like graph 
theory or within finite abstract simplicial complexes. We work here mostly in graph theory
and review the Morse inequalities $b(k)-b(k-1) + ... + b(0) \leq c(k)-c(k-1) +  ... + c(0)$
for the Betti numbers b(k) and the minimal number c(k) of Morse critical points of index k
and the Lusternik-Schnirelmann inequalities ${\rm cup}+1 \leq {\rm cat} \leq {\rm cri}$ 
between the algebraic cup length cup, the topological category cat and the analytic number
cri counting the minimal number of critical points of a function. 
\end{abstract}

\begin{center}
{\it In memory of Frank Josellis (1959-2022)}
\end{center} 


\section{Introduction}

\paragraph{}
Let $G=(V,E)$ be a {\bf finite simple graph} with {\bf vertex set} $V$ 
and {\bf edge set} $E$. A {\bf scalar function} $g$ on $G$ is a map 
$g: V \to \mathbb{R}$. 
The {\bf unit sphere} $S(v)$ of $v \in V$ is the sub graph generated by all 
vertices $w \in V$ with $(v,w) \in E$. 
The function $g$ is called {\bf locally injective} if $g(v) \neq g(w)$ for
all $w \in S(v)$. The {\bf stable unit sphere} $S^-_g(v)$ is the sub-graph of $S(v)$
defined by $S^-_g(v)=\{ w \in S(v), g(w)<g(v) \}$. A graph is {\bf contractible} 
if there is $v \in V$ such that both $S(v)$ and $G\setminus v$ are contractible, 
where $G \setminus v$ is the graph in which $v$ has been removed.
A vertex $v$ is a {\bf critical point} of $g$ if $S^-_g(v)$ is not contractible. 

\paragraph{}
The {\bf Lusternik-Schnirelman} category ${\rm cat}(G)$ of $G$ is the 
minimal number of contractible graphs that cover $G$. 
The minimal number of critical points which a locally injective
function can have is denoted by ${\rm cri}(G)$.
Graphs define a simplicial complex, an exterior derivative $d$, a 
{\bf Hodge Laplacian} $L=(d+d^*)^2$. The kernel of $L$ is the space of 
{\bf harmonic forms} which represents {\bf cohomology}. 
It carries a {\bf graded multiplication} called the {\rm cup product}.
The maximal number of positive degree forms which can be multiplied to get something non-zero
is the {\rm cup length} ${\rm cup}(G)$. Lusternik-Schnirelmann's theorem is
\fbox{${\rm cup}(G)+1 \leq {\rm cat}(G) \leq {\rm cri}(G)$}. 
\index{Graph ! finite simple}
\index{Finite simple graph}
\index{Function ! locally finite}
\index{Euler characteristic}
\index{simplex}
\index{face}
\index{clique}
\index{complete sub-graphs}
\index{simplicial complex}
\index{finite abstract simplicial complex}

\paragraph{}
A scalar function $g:V \to \mathbb{R}$ is {\bf Morse} if at every critical point $v$,
the graph $S_g^-(v)$ is a $k$-sphere for some $k$. The value $k+1$ 
is called the {\rm Morse index} of such a critical point $v$. If $v$ is a minimum of $g$ for 
example, then $S_g^-(v)=0$ is the empty graph $0$, which is a $-1$-sphere so that $k=0$. 
A {\bf $k$-sphere} $G$ is $k$-manifold for which $G \setminus v$
is contractible for some $v$. A {\bf $k$-manifold} is a graph for which every unit sphere 
$S(v)$ is a $(k-1)$-sphere. If $c_k(G)$ denotes the minimal number of critical 
points that a Morse function $g$ can have on $G$, and $b_k(G)$ is the $k$'th 
{\bf Betti number}, the dimension of the kernel of the Laplacian $L$ restricted to 
$k$-forms, then the {\bf strong Morse inequalities} are 
\fbox{$b_k(G)-b_{k-1}(G) + \cdots + (-1)^k b_0(G) \leq c_k(G)-c_{k-1}(G) +  \cdots + (-1)^k c_0(G)$}.
By adding the inequalities for $k$ and $k-1$, one gets the 
{\bf weak Morse inequalities} $b_k(G) \leq c_k(G)$. These inequalities work for all Morse functions
on a finite simple graph. 
\index{strong Morse inequalities}
\index{weak Morse inequalities}

\paragraph{}
Morse functions do exist on many graphs but not for all graphs.
The {\bf cube graph} or the {\bf dodecahedron graph}
for example are $1$-dimensional (triangle free) graphs. They both do 
not admit any Morse function. The reason is that for any $g$, the maximum has as 
$S_g^-(v)=S(v)= \overline{K_3}$ (the graph complement of the complete graph 
$K_3$ is a $0$-dimensional graph with 2 vertices and no edges),
which is not contractible so that the maximum is a critical point
but $S_g^-(v)$ is not a sphere. More generally, any vertex regular graph for which the 
unit sphere (which are all isomorphic graphs)
is not a sphere, can not carry a Morse function. 
\index{cube graph}
\index{dodecahedron graph}

\paragraph{}
Having seen this example, one could suspect that $G$ needs to be a manifold
in order to be able to carry a Morse function. This is not the case:
the graph $(V,E)$ defined by any finite abstract simplicial complex $G$ 
(in which $V=G$ and $E$ are the pairs $(x,y)$ with $x \subset y$ or $y \subset x$)
always has a Morse function. An example is the {\bf dimension function} 
$g(x) = {\rm dim}(x)=|x|-1$ which is locally injective because if $x$ is contained in $y$,
then $g(x) < g(y)$. It has the property that $S_g^-(x)$ is a sphere, the Barycentric refinement of
the boundary of the $k$-simplex $x$ is a $(k-1)$-sphere.
To see that $g$ is Morse, we build up the graph brick by brick. 
Start with the $0$-dimensional part, the vertices. These are all critical points of Morse index $0$
because the unit spheres are empty and so $(-1)$-dimensional spheres. 
Then we add edges which are critical points of Morse index $1$. Now add triangles, which are points 
of Morse index $2$ etc. 
The number $c_k(G)$ of critical points of index $k$ in $G$ agrees with $f_k(G)$, the number
of $k$-dimensional {\bf simplices} (=complete subgraph = {\bf clique}) 
$x$ in $G$. The vector $(f_0,f_1, \dots, f_d)$ is called the {\bf $f$-vector} of $G$. 
\index{f-vector}
\index{clique}
\index{complete subgraph}
\index{dimension functional}

\paragraph{}
We will look at Poincar\'e-Hopf and the Morse inequalities in the next section. But the just considered
example illustrates the theme. We have just seen that any simplicial complex has already a natural
Morse complex built in. Each cell is a ``handle". If $d$ is the maximal dimension of $G$, then the 
{\bf Poincar\'e-Hopf} formula $\sum_{k=0}^d (-1)^k c_k = \chi(G)$ 
rephrases with {\bf Euler-Poincar\'e} $\chi(G) = \sum_k (-1)^k b_k$ as the special case
$\sum_{k=0}^d (-1)^k b_k = \sum_{k=0}^d (-1)^k c_k$ of the Morse inequalities. The general 
Morse inequality appears by building up all simplices up to dimension
$k$. At each stage, {\bf Betti number} $b_k(G)$ is the {\bf nullity} of the 
{\bf $k$-form Laplacian} $L_k$ defined by 
$G$. If the simplicial complex of the graph $G$ has $n$ elements, then all exterior derivatives
are encoded in one lower triangular $n \times n$ matrix $d$ and $D=d+d^*$ is the {\bf Dirac matrix} 
of $G$. The {\bf Hodge Laplacian} $L= \oplus_{k=0}^d L_k$ is block diagonal and the block $L_k$
is the $k$-form Laplacian. 
\index{Euler Poincar\'e}
\index{dimension functional}
\index{Betti number}
\index{Dirac matrix}
\index{Form Laplacian}

\paragraph{}
The frame work of graphs is {\bf intuitive} and so also a paradigm for {\bf clarity}.
Every finite simple graph defines a finite geometry 
which allows to define discrete manifolds and so covers quite a bit of classical geometry. 
They are familiar because street networks, genealogy trees, social 
networks, computer networks or flow diagrams are graphs. 
We could start also with a {\bf finite abstract simplicial complex}, 
which stands for {\bf simplicity} as there is only one axiom:
it is a finite set of non-empty sets closed under the operation of taking finite non-empty subsets. 
Every graph defines such a complex. 
It is known as the {\bf Whitney complex}, {\bf flag complex} or {\bf order complex}. 
While we could start with a simplicial complex without
requiring it to come from a graph, it is more intuitive and so clear to use the language of graphs.
Restricting to graphs is almost no loss of generality because we can from a complex $G$
get a graph $(G,\mathcal{E})$, where $\mathcal{E}=\{ (x,y), x \subset y$
or $y \subset x \}$. There is a generalization of simplicial complexes called a {\bf delta set}. 
It is a paradigm for {\bf generality} as delta sets form a topos, a mathematical object that 
is closed under any operation like product, quotient, level sets. It includes many structures 
not covered by simplicial complexes like {\bf quivers}, graphs for which multiple-connections 
or loops are allowed. 
\index{quiver}
\index{delta set}
\index{simplicity, clarity, generality}
\index{Unix paradigm} 

\paragraph{}
Composing the {\bf Whitney map} with the {\bf graph formation map}
produces the {\bf Barycentric refinement}. According to Dieudonn\'e, the
Barycentric refinement is one of the three pillars of {\bf combinatorial topology};
the other two are {\bf incidence matrices} and {\bf duality}. 
Examples of complexes that are not Whitney complexes are
{\bf k-skeleton complexes} obtained by restricting a simplicial complex 
to simplices of dimension $\leq k$. The simplest example is
the $1$-skeleton complex of $K_3=\{ \{1,2,3\},\{1,2\},\{2,3\},\{1,3\},\{1\},\{2\},\{3\} \}$
which is $C_3=\{ \{1,2\},\{2,3\},\{1,3\}$, $\{1\},\{2\},\{3\} \}$. While no more
the Whitney complex of a graph, it is still a $1$-dimensional complex. Its Barycentric
refinement is $C_6$ which is now is the Whitney complex of the cyclic graph $C_6$.
Every complex which contains $C_3$ fails to be the Whitney complex of a graph. 
\index{Barycentric refinement}
\index{Whtney map}
\index{skeleton complex}

\paragraph{}
A simplicial complex $G = (G_0, \dots, G_n)$ defines {\bf face maps}
$d_{i,k}: G_k \to G_{k-1}$, 
$x \mapsto (x_0, \dots, \hat{x}_i,\dots,x_k)$ and so a {\bf delta set}, which is
an element of a {\bf finite topos}. 
A delta set is a presheaf over the simplex category $\Delta$ with strict inclusions as 
morphisms on $\Delta$. If all inclusions are chosen as morphisms, 
more conditions need to be satisfied for the functors and we get the smaller category
of {\bf simplicial sets}. Simplicial sets are more complicated in that there are two
type of structure maps which need more conditions. 
A {\bf delta set} is a geometry that is more general than simplicial complexes 
in that it is possible to do
products, build quotients, level sets of functions, where we can factor out symmetries or 
build Whitney type cell complexes from quivers. It even allows to model divisor structures
(geometric structures with additional integer valued structures telling about multiplicity). 

\paragraph{}
Every delta set has an exterior derivative $d$ and so has a cohomology, a Laplacian $L=(d+d^*)^2$,
leading to discrete partial differential equations like {\bf wave equation} $u_{tt}=-Lu$ or 
the {\bf heat equation} $u_t=-Lu$. 
The set $G=\bigcup G_k$ can be the vertex set of a graph,
where two elements are connected if there is a sequence of face maps
getting from one to the other. (The Dirac operator which contains the 
incidence matrices encodes this completely). This produces a finite simple graph which in general
has similar topological properties than $G$, at least if $G$ is a simplicial complex. 
As a data structure, we can encode a delta set as a triple $(G,D,r)$, where $G$ is a finite set
of $n$ elements, $D$ is a $n \times n$ matrix and $r:G \to \mathbb{N}$ is a dimension function. 

\paragraph{}
It is more convenient to ignore the face maps in the definition of a delta set and directly go to the 
exterior derivative $d$. A delta set is just a finite set of $n$ elements, a 
Dirac matrix $D = d+d^*$ and a dimension function which gives the dimension to 
each element $G_k={\rm dim}^{-1}(k)$ so that the elements of $G$ belonging to blocks of $L=D^2$ have 
fixed dimension. 
We need the dimension function $r$ because (unlike for simplicial complexes), for delta
sets, some of the dimension blocks $L_k$ of $L$ can be empty. The delta set $(G,D,r)=(\{1,2\},[0],(1))$
for example has only one element in $G$, the Dirac matrix is a $1 \times 1$ matrix and the element has
the assigned dimension $1$. No $0$-dimensional elements are present. 
The name ``Dirac matrix" has been chosen because it is a square root of the Hodge 
Laplacian $L$. It stands for a symmetric $n \times n$ matrix $D=d+d^*$ such that $d$ (and so 
the transpose) is nilpotent $d^2=0$. The trivial case $d=0$ reduces the delta set to a set. But $r$ 
can give a bit more structure. The open set $U=\{1,2\}$ in $K_2$ for example is a one-dimensional delta
set with Betti vector $(0,1)$. It is the complement of the closed set $K=\overline{K}_2$ in $G=K_2$ with
Betti vector $(2,0)$. Now $b(U) +b(K)=(2,0)+(0,1)=(2,1) > b(G)$ illustrates one of the simplest cases for
the fusion inequality. 
\index{face maps}
\index{cohomology}
\index{presheaf}
\index{topos}
\index{finite topos}
\index{wave equation}
\index{heat equation}
\index{partial differential equation}

\paragraph{}
A {\bf quiver} is a graph in which also multiple-connections and loops are allowed.
Quivers also carry an {\bf exterior derivative} $d$. One can use see them as $1$-dimensional
delta sets. If a quiver has $n$ nodes and $m$ edges, the {\bf gradient} $d$ maps a $0$-form to a $1$-form. It 
is a $m \times n$ matrix. A quiver is a delta set with face maps $d_{0,0},d_{0,1}: E \to V$.  
The category of quivers is a {\bf topos} too, as it is a presheave category on the Kronecker 
quiver $Q$. Quivers are natural as they can also be seen as a functor category. 
One can see them as $1$-dimensional delta sets but there are still $1$-dimensional delta sets 
which do not come from quivers.
Examples are $G=\{ \{1,2\} \}$ or  $G=\{ \{1\},\{1,2\} \}$ which both are open sets in the 
simplicial complex $K_2$. There is also a functor from {\bf quivers} to higher 
dimensional delta sets which extends
the Whitney functor from {\bf finite simple graphs} to finite abstract simplicial complexes.
\footnote{When talking about functors, one also has to take into account 
the morphisms. The standard morphisms in graph theory are graph 
homomorphisms, the usual morphisms in simplicial complexes are simplicial maps.}
Having multiple loops allowed in a graph was useful for spectral estimates, like in 
\cite{Eigenvaluebounds} because principal submatrices of quiver Laplacians 
remain quiver Laplacians. 
\index{topos}
\index{quiver}
\index{multigraph}

\paragraph{}
Delta sets are especially useful when looking at the {\bf topology} on a simplicial complex.
The open sets in a simplicial complex are the complements of sub-simplicial complexes. 
Open sets are no more simplicial complexes. But they are still delta sets. 
The cohomology of an open set is in general different from the cohomology of its
Barycentric refinement or closure: for $U=\{\{1,2,3\}\}$ for example, 
the {\bf Betti vector} is $(0,0,1)$, while the closure $G=K_3$ has Betti vector 
$(1,0,0)$. The complement $K$ 
of the open set $U$ is the skeleton complex $C_3$ which has Betti vector $(1,1,0)$ 
(it is a simplicial complex which is the $1$-skeleton complex of $K_3$ but not 
a Whitney complex of a graph). The relation $(1,1,1) = b(C_3)+b(U) \geq b(K)=(1,0,0)$ 
for Betti vectors is an example of the {\bf Fusion inequality} \cite{DiscreteAlgebraicSets}
which tells that $b(U) + b(K) \geq b(G)$ if $U,K$ are open-closed pairs in $G$, meaning
that $U$ is open, $K$ is closed and $U \cup K =G$ and $U \cap K = \emptyset$. 
The Barycentric refinement coming from the graph defined by $U$ is $K_1$ which has
cohomology $b=(1,0,0)$. 
\index{topology}
\index{Fusion inequality}
\index{1-skeleton complex}

\paragraph{}
In this review, which was written mostly in the fall of 2022 triggered by the passing of 
Frank Josellis (who had collaborated on this with me in \cite{josellisknill}, 
we keep it within graph theory and especially do not
attempt to generalize the  Lusternik-Schnirelmann theorem
nor the Morse inequalities to delta sets. Already when reformulating the results in the category 
of simplicial complexes, one has to adapt some definitions and make choices,
like whether to define the category as the minimal number of {\bf open sets} covering 
a simplicial complex or the minimal number of {\bf closed sets} covering it. 
The revision of the presentation benefited much from the work of Jennifer Gao 
\cite{Gao2024} who wrote in the spring of 2024 a senior thesis about topological 
aspects of graph theory.

\paragraph{}
Note that there are various different definitions of Lusternik-Schnirelmann category in the 
literature in the continuum. There are different ways one can proceed to carry them over to the
discrete. When working on \cite{josellisknill}, we started with other 
definitions than what is used here and tried to get to notions that are homotopy invariant. 
We eventually decided not to worry about homotopy invariance as the notion of category is a rather
topological notion. A difficulty is that category as defined here is not homotopy invariant. 
The {\bf dunce hat} is not contractible so that the category is larger than 1. But it is 
homotopic to $1$. There must therefore have been a homotopy step which changed the Lusternik-Schnirelmann
category as defined here. 

\paragraph{}
When defining Lusternik-Schnirelmann category one has to decide first of all, what  elements should be used 
as coverings. One can use open contractible sets, closed contractible sets, open sets that have closures
homotopic to 1 or closed sets that are homotopic to 1. The decision to take contractible make the notion 
computable. One could also look at the notion of {\bf contractible within an other host space}: 
an example is a closed path in a graph. It is intrinsically never contractible
as the cohomology of a circle is $b=(1,1)$ and not zero; but a closed loop in a contractible space
is considered to be {\bf homotopically trivial} like when defining the {\bf fundamental group}. 
For us, {\bf contractible parts} of a cover always means intrinsically contractible. 
\index{Frank Josellis}

\paragraph{}
We had attempted to summarize some aspects of geometry and calculus in graphs in 
expository documents \cite{knillcalculus,KnillBaltimore,AmazingWorld}. 
The definitions we have used in the last
decade have fluctuated a bit, even for basic things like the notion of a {\bf manifold}, 
the definition of {\bf Morse function} or the choice of definition for 
Lusternik-Schnirelmann category. The current choice of definitions has remained
pretty stable. It is guided by three ``functionals":
{\bf the simplicity of the definition} and that theorems should work 
{\bf as close as possible to the continuum} without much modification. 
Finally, the {\bf proofs should be simple}. 

\paragraph{}
We realized over the years that many graph theorists think about
graphs differently. A common view is to see them as 
{\bf $1$-dimensional simplicial complexes}. It should be clear
however, that the language of graphs is just one of many approaches to cover
material which usually is described using the language of simplicial complexes or
notions from the continuum. Already in older texts like \cite{AlexandroffHopf},
one can see that graphs were drawn to visualize continuum manifolds. In computer
graphics, surfaces are drawn as triangulations and the graphics complexes only 
need vertices and edges as well as the position of the vertices in space.

\paragraph{}
In many texts that use simplicial complexes, also describes "usual simplicial complexes"
meaning that one looks at notions in the continuum.
Or then one uses the {\bf geometric realization} functor giving from an abstract version
a Euclidean realization. Even Dehn and Heegaard who introduced
abstract simplicial complexes, use also the continuum \cite{DehnHeegaard}. 
In this text, we never look at geometric realizations. The continuum 
is not wanted. We want to avoid the infinity axiom. 
Everything is combinatorial. While we mention the real line 
$\mathbb{R}$ like when talking about $\mathbb{R}$-valued functions,
it should be clear that one could very well also use integer-valued 
functions or functions taking values in a finite totally ordered set. 
We always can remain in a finitist setting. 
\index{geometric realization}
\index{finitist}
\index{infinity axiom}

\section{Poincar\'e-Hopf}

\paragraph{}
The {\bf Euler characteristic} of a graph $G$ is defined as
the Euler characteristic of its Whitney complex:
$\sum_{x \in G} \omega(x)$ with $\omega(x)=(-1)^{{\rm dim}(x)}$.
The 1-dimensional Euler characteristic $|V|-|E|$ would just take 
into account the $1$-dimensional skeleton complex. The historically
best known situation is the 2-dimensional $|V|-|E|+|F|$, where $F$ 
is the set of 2-dimensional faces. This formula has been discovered 
by Descartes \cite{Aczel}. For convex polyhedra (to which the majority
of books on polyhedra restrict to), there is no difficulty 
to see the highest dimensional facets. Especially in the case of 
2-dimensions, where we deal with planar graphs. 

\paragraph{}
The simplicial complex for the graph $(V,E)$ has as sets the vertex sets of 
complete sub-graphs of the graph. We denote it again $G$. For a general delta 
set $G$, the Euler characteristic is defined in the same way. It is 
$\chi(G) = \sum_{j=1}^n (-1)^j |G_j|$, where $G_j$ 
are the elements of dimension $j$. 

\paragraph{}
Given a locally injective function $g$ on the vertex set $V$ of a graph,
define $S^-_g(v)$ as the sub-graph generated by all vertices in the sub-graph 
$S(v)$, where $g(w)$ is smaller than $g(v)$. 
One can now define {\bf Poincar\'e Hopf indices} $i_g(v) = 1-\chi(S^-_g(v))$.
As an integer-valued function on the set $V$ of vertices of $G$, it is a {\bf divisor}
in the sense of the {\bf discrete Riemann-Roch theory} \cite{BakerNorine2007}. 
\index{Euler characteristic}
\index{divisor}
\index{Riemann-Roch}
\index{discrete Riemann-Roch}

\paragraph{}
Note that we look here at functions on the {\bf vertex set} $V$ of the graph $G$ and 
not at functions on the simplicial complex $G$, which is the vertex set of 
the Barycentric refinement graph $G_1$. An example for the later is the
function $\omega(x)=(-1)^{{\rm dim}(x)}$. It can be seen as the {\bf Poincar\'e-Hopf index function} of 
the dimension function $g(x)={\rm dim}(x)$. The later is a locally injective function 
on the Barycentric refinement $G_1$. The sum over all these Poincar\'e-Hopf indices
on $G_1$ is the definition of the Euler characteristic of $G$. 
The statement that the total Poincar\'e-Hopf indices of $g$ adds up to the Euler 
characteristic is already an example of the Poincar\'e-Hopf theorem. 

\paragraph{}
By induction in the number of vertices, one immediately gets the {\bf Poincar\'e-Hopf formula}. 
It is a discrete analog of the continuum \cite{poincare85,hopf26,Spivak1999}, in the special case
if the vector field is a gradient field. We will see below how to generalize this to ``vector fields"
like a digraph structure on the graph having the property that there is no closed cycle 
in each simplex. We formulated the following result in \cite{poincarehopf} with a more complicated
proof. It matured while working on Lusternik-Schnirelmann in \cite{josellisknill}.

\begin{thm}[Poincar\'e-Hopf]
If $g$ locally injective then $\chi(G) = \sum_{v \in V} i_g(v)$ 
\label{PoincareHopf}
\end{thm}
\begin{proof}
The induction foundation is the one-point graph $1=K_1$, with $\chi(K_1)=1$, 
where $i_g(v)=1$ because the unit sphere $S^-_g(v)$ is empty independent of the function $g$. 
Now to the induction step:
if the result is known for all graphs with $(n-1)$ vertices and all colorings 
(locally injective functions on $V$) on such graphs,
and $G$ is a graph with $n$ vertices with a coloring $g$, identify a local maximum $v$
of $g$. The graph $G \setminus v$ has $n-1$ elements and induction assumption can be used. 
Then $\chi(G)=\chi(G \setminus v) + \chi(B(v))-\chi(S(v))= [\sum_{w \neq v} i_g(w)] + i_g(v)$, 
because every unit ball $B(v)$ has $\chi(B(v))=1$ and $S_g^-(v) = S(v)$ if $v$ is a 
local maximum of $g$. Note that adding the new vertex $v$ has not altered 
the graphs $S_g^-(w)$ of the other vertices so that the indices $i_g(w)$ stayed the same. 
This proves the induction step. 
\end{proof} 

\index{Poincar\'e-Hopf theorem}
\index{Theorem ! Poincar\'e-Hopf}

\paragraph{}
There are various way how one can generalize this result 
\cite{poincarehopf,MorePoincareHopf,PoincareHopfVectorFields,parametrizedpoincarehopf}. 
Instead of the Euler characteristic,
one can look at the {\bf $f$-function} $f_G(t) = 1+f_0 t + f_1 t^2 + \cdots f_d t^{d+1}$
which satisfies $\chi(G)=1-f_G(-1)$. The statement is 
$f_G(t) = 1+t \sum_{x \in V} f_{S_g^-(x)}(t)$, where $S_g^-(x) = \{ y \in S(x), g(y)<g(x) \}$
One can also replace the Euler characteristic by an {\bf energy sum} 
$H(G) = \sum_{x \in G} h(x)$, if $h$ is a ring-valued function on the simplicial complex. 
This generalizes the situation when $h(x)=\omega(x)$.
\index{energy sum} 
\index{f-function}

\paragraph{}
Savana Ammons \cite{SavanaAmmons} informed us in 2023 about a previously overlooked
Poincar\'e-Hopf result for graphs on classical 2-manifolds \cite{Glass1973}. 
It can be generalized to discrete 2-dimensional CW complexes (which can be seen as delta sets)
and which does not involve the continuum. 
In two dimensions, the traditional view is to see $G$ embedded in a $2$-dimensional surface. 
One has then a notion of ``face". We note however that Theorem~\ref{PoincareHopf} works in arbitrary dimensions and is 
close to the traditional Poincar\'e-Hopf theorem for gradient fields. We have formulated
\cite{MorePoincareHopf} also with digraphs without circular parts in any simplex as this produces
a total order on each simplex and so a {\bf vector field} $F: G \to V$, where $G$ is the simplicial 
complex and $V$ the vertex set. Directing every energy $h(x)$ contributing to the Euler characteristic
$\chi(G) = \sum_{x \in G} h(x)$ to its vertex produces a function $i_F(v)$ on vertices which satisfies 
$\sum_v i_F(v)=\chi(G)$. 

\paragraph{}
The just given argument works also for the above proven theorem, if $g$ is a
locally injective function on the vertices. Being locally injective
makes sure that $g$ defines a {\bf total order} on each simplex $x \in G$. 
We can now move every energy value $\omega(x)$ entering 
the Euler characteristic $\sum_{x \in G} \omega(x)$ to the
vertex in $x$, where $g$ is maximal. The sum of the indices is the same;
indeed, no $\omega(x)$ has been lost nor double counted.

\paragraph{}
The frame work can be extended a bit more. 
We can think of an assignment which maps a simplex $x$ to a vertex $v = F(x) \in x$ as
a {\bf vector field}. It is a generalization of a {\bf direction} because on a 1-dimensional 
simplicial complex, we can see $F$ as defining an orientation of the edges. 
For $1$-dimensional simplicial complexes, meaning a graph with the $1$-dimensional 
skeleton complex a vector field is the same than assigning a direction to each edge.
The vector field structure so upgrades the graph to a {\bf directed graph}. 
\index{direction}
\index{vector field}
\index{total order}

\paragraph{}
The just given argument for Poincar\'e-Hopf is in the spirit of Morse theory:
pick the vertex $v$ on which the function $g$ is maximal, 
use that $\chi(B(v))=1$ for any {\bf unit ball} $B(v)$ of $G$ and the {\bf valuation formula} 
$\chi(A \cup B) = \chi(A) + \chi(B) - \chi(A \cap B)$ which is valid for simplicial 
complexes as well as for sub-graphs $A,B$ of $G$ that generate themselves in $G$.
An injective function $g$ defines a {\bf Morse type build up} of the graph:
we have a set of graphs $c \to G(c)$ generated by vertices $V(c)=\{ v, g(v) \leq c \}$. 
Points $v$, where $\chi(S_g^-(v)) \neq 1$, are automatically critical points because
the Euler characteristic changes by $1-\chi(S_g^-(v))$. If $S_g^-(v)$ are $k$-spheres, 
meaning that $f$ is a {\bf Morse function}, the changes are by $(-1)^{m(v)}$, 
where $m(v)-1$ is the dimension of the sphere $S_f^-$. 
For a minimum $v$ of $f$ for example, $S^-_f(v)=0$, the $(-1)$-dimensional sphere
so that the {\bf Morse index} $m(v)=0$ and the Poincar\'e-Hopf index $i_g(v) = 0$.
\index{Morse index}
\index{Poincar\'e-Hopf index}
\index{Valuation formula}

\paragraph{}
Let us look a bit more at the {\bf Glass theorem} \cite{Glass1973}, which is the special
case of a $2$-dimensional complex but where one has also faces coming from an embedding 
in a surface. We can see the embedded graph structure as a delta set
$(V,E,F)$, where $F$ are the {\bf faces}, connected regions enclosed by edges. If we
prefer to stay in combinatorics and avoid the continuum, we could just add the faces
and require that the Barycentric refinement of this structure is a discrete 2-manifold, 
a graph for which every unit sphere is a circular graph with 4 or more elements. 
A given (fixed but arbitrary) orientation on each of the unit spheres $S(v)$ now 
allows to redistribute the values $\omega(x)$ for
each of the simplices $V \cup E \cup F$.  Keep the values $\omega(x)$ from $V$
and $F$ where they are then move the value $\omega(e)$ of an edge $e=(a,b)$ 
either to the vertex set $b$ or then to the face $f=(a=a_0,b=a_1,\dots,a_{n-1})$ 
defined by the orientation. For a vertex $x=v$ or face $x=v$, the {\bf Glass index} is 
$1-R(x)/2$, where $R(x)$ is the number of orientation changes of edges entering 
or surrounding $x$. The division by $2$ comes from the fact that every edge $e$ induces
exactly two sign changes on neighboring part of $V \cup F$. 

\index{Poincar\'e-Hopf theorem}
\index{Glass theorem}
\index{unit ball}
\index{valuation} 
\index{Morse buildup}
\index{Morse index}

\paragraph{}
Let us reproof quickly the Poincar\'e-Hopf version for {\bf digraphs} (directed graph)
$G=(V,E)$ as given in \cite{PoincareHopfVectorFields}.
Let us call a {\bf digraph} {\bf locally non-circular} if no closed path exists in 
any any of its cliques $x$ (=simplex=face=complete subgraph). 
This especially means that if $(a,b) \in E$, then 
$(b,a)$ is not in $E$. The relation $E$ now defines a {\bf total order} on each 
simplex. Now think $\omega(x)=(-1)^{{\rm dim}(x)}$ as an {\bf energy} attached to the simplex $x$. 
The Euler characteristic is the total energy. One can now in any simplex identify the 
maximal vertex $v$ given by the total order and assign the energy $\omega(x)$ to $v$. If
this is done for all simplices, we get a function $i(v)$ on the vertex set and naturally
$\sum_{v \in V} i(v)$ is equal to $\sum_x \omega(x)=\chi(G)$.  We have shown: 
\index{digraph}
\index{locally non-circular}
\index{non-circular}
\index{Poincar\'e-Hopf digraph}

\begin{thm}[Poincar\'e-Hopf for digraphs]
For a locally non-circular digraph, $\sum_{v \in V} i(v) = \chi(G)$. 
\end{thm}

\paragraph{}
The paper \cite{PoincareHopfVectorFields} generalized this to more abstract vector fields
on a graph $(V,E)$ with simplex set $G$. 
A {\bf vector field} is just a map $F: G \to V$ with $F(x) \in x$. One can enhance this 
with an additional map $i: V \to G$, where $v \in i(v)$ and so get a {\bf permutation}
$T=i \circ F: G \to G$. A vector field then not only defines the Poincar\'e-Hopf indices
but also a {\bf dynamics}. Equilibria are simplices $x$ for which $T(x)=i \circ F(x)=x$. 
The permutations defined by the vector field $F$ model now more closely a vector field
dynamics in the continuum because $T(x)$ intersects with $x$. 
\index{vector field} 
\index{dynamical systems} 
\index{equilibrium} 

\section{Gauss-Bonnet}

\paragraph{}
Poincar\'e-Hopf indices are {\bf divisors} on the graph (maps $V \to \mathbb{Z}$).
We will see curvature as an average of such indices. It is
in general not a divisor because curvature values are in general
rational. {\bf Integral geometry} links Poincar\'e-Hopf indices with curvature.
It uses probability theory to define geometric quantities
like length, volume or curvature. It is a classical topic \cite{Weyl1939,blaschke,
Santalo,Gray, Banchoff67,Banchoff70}.
Given any {\bf probability space} $(\Omega,\mathcal{A},m)$ 
on the set of {\bf locally injective functions} $\Omega$, one gets a
{\bf curvature} $K(x)={\rm E}[i_g]$ as an {\bf index expectation} 
\cite{indexexpectation,indexformula,DiscreteHopf2,
ConstantExpectationCurvature,DiscreteHopf,eveneuler}.
Because of Fubini's theorem, the expectation commutes with summation and
curvature automatically satisfies the {\bf Gauss-Bonnet formula}. 
``Gauss-Bonnet is an expectation of Poincar\'e-Hopf". 
\index{probability space}
\index{Fubini theorem}
\index{index expectation}
\index{integral geometry}

\begin{thm}[Gauss=Bonnet]
For any index expectation curvature, we have  $\chi(G)=\sum_{v \in V} K(v)$. 
\end{thm}

\begin{proof}
This directly follows from Fubini: we can switch the finite
summation and expectation procedures: 
$\sum_v K(v) = \sum_{v \in V} {\rm E}[i(v)] = {\rm E}[ \sum_{v \in V} i(v) ] 
= {\rm E}[\chi(G)] = \chi(G)$. 
\end{proof} 

\paragraph{}
We stepped into this topic first with \cite{elemente11,cherngaussbonnet,knillcalculus}. 
For graphs which define 2-dimensional manifolds, the curvature is $K(v)=1-deg(v)/6$
which is a century old theme, especially in the graph coloring literature
(see \cite{knillgraphcoloring}) and going back to \cite{Eberhard1891}. 
If the measure is the counting measure on all locally injective functions for 
example, one gets the {\bf Levitt curvature} 
$K(v) = 1-f_0(S(v))/2+f_1(S(v))/3 - \dots$, where $S(v)$ is the 
{\bf unit sphere} of $v$, the graph generated by all vertices attached to $v$ 
\cite{Levitt1992}. For 2-dimensional manifolds, where 
$f_0(S(v))=f_1(S(v))={\rm deg}(v)$, we have $K(v)=1-{\rm deg}(v)/6$. For
1-dimensional graphs, graphs without triangles, we have 
$K(v)=1-{\rm deg}(v)/2$. 
\footnote{We noticed \cite{Levitt1992} only after \cite{cherngaussbonnet} was written. 
\cite{Levitt1992} does not stress the Gauss-Bonnet connection.}

\paragraph{}
The Poincar\'e-Hopf theorem can be understood as a rule which pushes the signed measure 
(= energy) $\omega(x) = (-1)^{{\rm dim}(x)}$ from the simplix $x$ to 
a vertex $v$ in $x$, where $f$ is largest ending up with a weight $i_g(v)$ 
on the vertex $x$. The symmetric Gauss-Bonnet formula of Levitt can then be
understood as distributing $\mu(x)$ equally to all $(k+1)$ 
vertices that are present in the $k$-simplex $x$. 
All this works in more generality, when $\omega(x)$ can be a quite arbitrary function 
to a ring and especially for functions that take values in the units of real 
division algebras. 
\cite{EnergizedSimplicialComplexes2,EnergizedSimplicialComplexes}. 

\index{Gauss-Bonnet formula}
\index{curvature}
\index{division algebra valued functions}
\index{integral geometry}
\index{Levitt formula}

\paragraph{}
Having a probability measure on locally injective functions is a way to {\bf deform} the 
geometry. We think about it as a way to impose a {\bf discrete Riemannian metric}. 
The reason is that a probability space on functions also defines a distance function 
on $V$. Let $v,w$ be two vertices in a graph, its {\bf Crofton distance} 
is the probability of the set of
functions $g$ for which $g(v)$ and $g(w)$ have different signs. This is part of integral geometry 
\cite{Santalo,Santalo1,KlainRota,Schneider1}. The idea of index expectation in the context of 
Gauss-Bonnet probably first appeared in \cite{Banchoff67,Banchoff70} even-so one can see the ideas
of integral geometry in the Gauss-Bonnet-Chern theme in general.
In the discrete, integral geometry paves a way to add a metric and curvature structure which 
is {\bf deformable}. Similarly as one can deform the Riemannian metric in the continuum, we can tune
the probabilities on functions. This allows more flexibility and hope that using suitable notions of
``sectional curvature" in discrete manifolds, one could get a Gauss-Bonnet curvature which 
is positive for positive curvature manifold. See 
\cite{ConstantExpectationCurvature,DiscreteHopf2,DiscreteHopf}.
\index{Crofton formula}
\index{Integral geometry}
\index{discrete Riemannian metric}

\paragraph{}
The general {\bf integral geometric definition of distance} on a Riemannian manifold goes 
as follows: if $\gamma: [a,b] \to M, t \to r(t)$ is a smooth curve in a smooth connected manifold $M$,
denote by $L(\gamma)=\int_a^b |r'(t)| \; dt$ its length (where $ds=|r'(t)| dt$ is measured using 
the Riemannian metric). Given a probability space
$(\Omega,\mathcal,P)$ of smooth Morse functions $\omega$ on $M$, we can look at the random variable
$N_{\gamma}(\omega)$ counting the number of intersections
of the level surface $\{ \omega=0 \}$ with $\gamma$. More precisely, we count the number of transitions
from $f \leq 0$ to $f>0$. This defines a {\bf pseudo metric}
$$ d(x,y) = \inf_{\gamma(x,y), N_\gamma \in L^1(\Omega,P)} {\rm E}[N_{\gamma}]  \; , $$
where the infimum is taken over all curves connecting $x$ with $y$ with the understanding that
$d(x,y) = \infty$, if there should be no $\gamma$ for which $N_\gamma$ is in $L^1$.
The {\bf Kolmogorov quotient} $(M_P,d_P)$
consists of all equivalence classes, with the equivalence relation given by $x \sim y$ if 
$d(x,y)=0$. This is a generalization of Riemannian metric because if we take a Riemannian manifold
and isometrically embed it into an ambient space \cite{EssentialNash}
and use a rotationally symmetric measure, then we
recover the standard metric. 
\index{Kolmogorov quotient}
\index{pseudo metric}
\index{generalized Riemannian metric}

\section{Homotopy} 

\paragraph{}
The {\bf unit sphere} $S(v)$ of a vertex $v$ is the sub-graph of $G$ 
generated by all $w \in V$ with $(v,w) \in E$. In the language of {\bf metric
spaces}, it is the unit sphere with respect to the {\bf geodesic metric} on $G$,
where adjacent nodes have distance $1$. But in graph theory, it is understood 
as a graph and not just as a set of vertices. 

\paragraph{}
A graph $G=(V,E)$ is called {\bf contractible} if there exists $v \in V$,
such that both the unit sphere $S(v)$ and the graph $G \setminus v$, 
the sub-graph of $G$ generated by all vertices of $G$ different from $v$, 
are contractible. This inductive definition is primed by the 
assumption that the {\bf one-point graph} $K_1$ is contractible. The empty
graph is not considered contractible. Note that the definition also depends on 
the simplicial complex we have chosen. If we would look at graphs as one-dimensional 
simplicial complexes, then already a triangle $K_3$ would not be contractible. 

\paragraph{}
What is here called contractible is often also called collapsible. 
If we can use homotopies (=contraction and extension steps) 
to get from a graph to $K_1$, we call this {\bf homotopic to $1$}. 
\footnote{Rather than using the ``gotcha" nomenclature of contractible and collapsible
which is used inconsistently in the literature anyway, we chose to identify 
contractible and collapsible and use "homotopic to 1" for the wider equivalence
relation. It makes sense to stress this point and repeat it. The difference
between contractible and "homotopic to 1" could not be bigger, justifying to distinguish
it also with nomenclature.}
\index{contractible}
\index{collapsible}
\index{homotopic to 1}
\index{geodesic metric}
\index{unit sphere}

\paragraph{}
This definitions of manifold was chosen so that the notion of ``sphere" and so the 
notion of ``manifold" is {\bf constructive}. One could replace ``contractible" with 
``homotopic to $1$", but that would be harder to check and lead to a definition which 
can not be checked with a fixed algorithm The inductive definition of
sphere goes back to Evako, formerly Ivashchenko \cite{Ivashchenko1993,Evako1994,I94a}.
The process of removing a vertex $v$ with contractible $S(x)$ or the reverse
process of linking a new point to a contractible part of a graph
are what constitutes {\bf homotopy steps}. Two graphs are called {\bf homotopic}, 
if a sequence of homotopy steps allows to get from one to an other. 
Unlike checking whether a graph is contractible, checking whether a graph is 
{\bf homotopic to $1=K_1$} can be hard. The {\bf dunce hat graph} is an example 
of a graph which is not contractible but which is {\bf homotopic to $1$}: 
it is a graph which needs to be expanded first before it can be contracted. For more about
the difficulty of deciding whether a structure is a sphere or not, see
\cite{Novikov1955,Chernavsky,SphereRecognition}. 
\index{constructive}
\index{dunce hat}
\index{decision problem}
\index{homotopy step}
\index{dunce hat}

\paragraph{}
Homotopy in discrete finite frame work has a long history. 
Historically,  combinatorial notions of homotopy were put forward by 
J.H.C. Whitehead \cite{Whitehead}. For graphs, the process of homotopy has been defined in \cite{I94} 
and was refined in \cite{CYY} (we consider the later a significant step as it 
simplified the notion considerably). These are definitions which work when we look at
the Whitney complex. Homotopy deformations of graphs produce deformations of the
simplicial complexes which would produce homotopy deformations of geometric
realizations. We can see homotopy deformations also as the process of adding a disjoint open 
set $U$ of trivial cohomology $\vec{b}(U)=\vec{0}$ to a closed set
\cite{DiscreteAlgebraicSets}: $G \to G \cup U$ produces again a simplicial complex
which by the {\bf trivial cohomology} of $U$ must have the same cohomology. 
The space $U=\{ \{1\},\{1,2\} \}$ for example has trivial cohomology 
(something which is not possible for simplicial complexes). Attaching this
to the leaf of a one dimensional graph produces a homotopy extension of the graph.
See \cite{May2008,MayPishevar2008} or finite topological spaces. 
\index{unit sphere}
\index{contractible graph}
\index{one-point graph}
\index{homotopy step}
\index{dunce hat}

\paragraph{}
To summarize, unlike the straightforward {\bf contractibility}, 
{\bf homotopy} is an equivalence relation on the category of 
finite simple graphs that is computationally difficult. 
To decide whether two graphs are homotopic
or not might, can need some creativity. To see for example that the graphs 
$C_5$ and $C_6$ are homotopic, we first have to ``fatten" it. It needs a few 
homotopy steps. But like all discrete manifolds (a notion covered in the next
section), these graphs $C_5$ and $C_6$ are not contractible. It is impossible to 
contract $C_6$ to $C_5$ but it is possible to make a homotopy deformation 
from $C_6$ to $C_5$. In a general manifold, one can not even remove one single vertex,
as every unit sphere is a sphere
and is not contractible. Already the $(-1)$ sphere $0$ and 
the {\bf 0-sphere}, the 2-point graph $\overline{K_2}$, the graph complement of the 
complete graph $K_2)$, are not contractible. 

\section{Manifolds}  

\paragraph{}
A graph $G=(V,E)$ is a {\bf discrete $d$-manifold} if every {\bf unit sphere}
$S(v)$ is a $(d-1)$-sphere. Inductively, a {\bf $d$-sphere} is a discrete 
$n$-manifold if there exists a vertex $v$ such that $G \setminus v$ is 
{\bf contractible}. The empty graph $0$ is declared to be the $(-1)$-sphere. 
So, the $2$-vertex graph $\overline{K_2}$ without any edge is the {\bf $0$-sphere} and
a cyclic graph with $4$ or more vertices is a {\bf $1$-sphere}.
A graph is called a {\bf $d$-ball} if it is of the form $G \setminus v$, where $G$ is 
a $d$-sphere. A graph is a {\bf d-manifold with boundary}, if {\bf (i)} every point 
has a unit sphere that is either a $(d-1)$-sphere or a $(d-1)$-ball and 
{\bf (ii)} both cases do appear. A complete graph $K_{d+1}$ for example is not a 
$d$-manifold with boundary but its Barycentric refinement is an $d$-ball.
The refinement of the $2$-simplex $K_3$ for example is a wheel graph 
with $7$ vertices. The simplicial complex of $K_3$ is 
$\{ \{1\},\{2\},\{3\},\{1,2\},\{2,3\},\{1,3\},\{1,2,3\}\}$, a set with seven elements. 

\paragraph{}
The {\bf boundary operation} $G \to \delta G$ can be defined for a general graph $G$.
The graph $\delta G$ is the sub-graph generated by the set of vertices for which the unit 
sphere $S(v)$ is contractible.  The remaining vertices are called {\bf interior points}.
A graph is called a {\bf d-manifold with boundary} if every interior point 
$v$ has a unit sphere $S(v)$ that is a $(d-1)$-sphere and every boundary point 
$v$ has a unit sphere that is a {\bf $(d-1)$-ball}, (which by definition is 
a {\bf punctured $(d-1)$ sphere}). An example of a 2-manifold with boundary is a 
wheel graph $W_n$ with boundary $C_n$.
\index{boundary ! general graph}
\index{puncture sphere}
\index{ball} 

\paragraph{}
The complement of the boundary in $G$ is the {\bf interior} of the manifold. It
it the graph generated by the interior points. There is a more sophisticated notion 
of {\bf interior} in a simplicial complex by taking the union of all stars of 
vertices away from the boundary and where the boundary is the set of simplices
for which $S(x) = \overline{U(x)} \setminus U(x)$ with star $U(x)$ is a contractible
subcomplex. 

\paragraph{}
Here, in a graph theoretical setting, the {\bf interior} is the set of 
points $v$ for which the unit sphere $S(v)$ is not contractible. In a manifold
with manifold setting, it is the set of $v$ for which $S(v)$ is a sphere.  
We always assume that a manifold with boundary both has a non-empty interior 
and a non-empty boundary. 
The unit sphere $S_{\delta G}(v)$ of every $v \in \delta G$ is the 
boundary of the ball $S_G(v)$ and so a $(d-2)$-sphere. 
Therefore, if $G$ is a $d$-manifold with boundary $\delta G$, 
then its boundary $\delta G$ is a $(d-1)$-manifold without boundary.
\index{discrete manifold}
\index{boundary}
\index{interior}
\index{manifold}
\index{manifold with boundary}
\index{(-1)-sphere}

\paragraph{}
A $d$-{\bf ball} $B$ is defined to be a graph which is of the form 
$S \setminus v$ for some $d$-sphere $S$ and vertex $v$, where $S \setminus v$
is the subgraph generated by all vertices different from $v$.
Every ball $B$ is contractible because by definition, 
$B=S \setminus v$ for {\bf some} $v$ is contractible. 
If $S$ is a $d$-sphere, then also $B=S \setminus w$ for any other vertex $w$ is 
contractible: to see this, we only have to show that if $S \setminus v$ is contractible
and $w$ is a neighboring vertex, then $S \setminus w$ is contractible. 
For any ball $B$, we can by definition 
remove a boundary vertex, unless it is the {\bf unit ball} 
$B(v)=S(v) \cup \{v\}$ of a single vertex $v$. The unit ball $B(v)$ is
contractible but removing any vertex on the boundary produces a graph 
without interior. 

\paragraph{}
We can check by induction that every contractible graph $A$ has 
Euler characteristic $\chi(A)=1$ and 
that every $d$-sphere $A$ satisfy the {\bf Euler gem formula} 
$\chi(A) = 1+(-1)^d$. See \cite{lakatos,Richeson}.
The induction step uses the {\bf valuation formula} 
$\chi(A \cup B) =\chi(A)+\chi(B)-\chi(A \cap B)$ which holds
for two arbitrary sub-graphs $A,B$ that generate themselves in $G$. 
Here is the induction step: if $A$ is a $d$-sphere, then $S(v)$ is a 
$(d-1)$-sphere that generates itself and $B(v)$ is a $d$-ball that generates
itself. Therefore $\chi(A) = \chi(A \setminus v) + \chi(B(v)) - \chi(S(v)) 
         = 1+1-(1+(-1)^{d-1})=1+(-1)^d$.
\index{Euler gem formula}
\index{valuation formula}

\paragraph{}
A discrete $d$-manifold is never contractible because no vertex can be removed: 
indeed, every unit sphere $S(v)$ is a $(d-1)$-sphere and since $\chi(S(v)) \in \{0,2\}$, 
we would remove a point with non-zero index $1$ or $-1$. There is no homotopy reduction
possible because a homotopy reduction step would remove 
a point with contractible unit sphere. But the unit sphere of a $d$-manifold is a
$(d-1)$-sphere and so not contractible. That a d-sphere is not contractible follows
from the just computed Euler gem formula for Euler characteristic 
(either $\{0,2\}$) and that a contractible graph has Euler characteristic $1$. 

\index{unit ball}
\index{Euler gem formula}
\index{valuation formula}
\index{gem formula}
\index{Euler characteristic of a sphere}

\section{Dimension}

\paragraph{}
The {\bf maximal dimension} ${\rm dim}_{max}(G)$ of a graph $(G,V)$ is 
$k$, if the graph contains a complete sub-graph $x=K_{k+1}$ but does not contain
any larger complete sub-graph $K_{k+2}$. 
Complete sub-graphs $x$ are also called {\bf faces}, {\bf cells},
{\bf simplices} or {\bf cliques} (where the later terminology is mostly used 
in graph theory). The set of all vertex sets of complete 
sub-graphs $x$ in $G$ is a {\bf finite abstract simplicial complex} $G$,
a finite set of sets $x$ closed under the operation of taking 
finite non-empty subsets. The integer ${\rm dim}_{max}(G)+1$
is also known as the {\bf clique number} of $G$. 
The maximal simplices are also known as {\bf facets}. 
The maximal dimension of the {\bf empty graph} $0$ is assumed to be $(-1)$. The 
empty graph is the {\bf initial object} in the category of graphs; it is an object of
Euler characteristic $0$. The $1$-point graph $1$ is the {\bf terminal object} in the 
category. 
\index{maximal dimension}
\index{clique number}
\index{empty graph}
\index{initial object}
\index{terminal object}
\index{facet} 

\paragraph{}
The {\bf inductive dimension} of a vertex $v \in V$
in a graph $G$ is defined as ${\rm dim}_G(v)=1+{\rm dim}_{ind}(S(v))$,
where ${\rm dim}_{\rm ind}(A)$ is the average of all inductive 
dimensions of vertices in $A$. With this definition,
${\rm dim}_{ind}(G) = \frac{1}{|V|} \sum_{v \in V} {\rm dim}_{ind}(S(v))$,
one has ${\rm dim}_{ind}(A) \leq {\rm dim}_{max}(A)$. 
Examples, where equality holds are discrete manifolds 
with or without boundary, regular graphs and especially complete graphs. 
The inductive dimension behaves a bit like the Hausdorff dimension
for metric spaces. The inductive dimension is bounded above by the maximal
dimension but the discrepancy can be arbitrarily large as one can see by 
adding as many zero dimensional parts to a given graph (which lowers the average
dimension but does not change the maximal dimension.) 
\index{inductive dimension}
\index{dimension ! inductive}
\index{maximal dimension}
\index{dimension ! maximal}

\paragraph{}
The {\bf augmented dimension} is defined as the inductive dimension incremented by $1$:
${\rm dim}^+(G) = 1 + {\rm dim}(G)$. In \cite{BetreSalinger} it was shown that
that ${\rm dim}^+(G+H) = {\rm dim}^+(G) + {\rm dim}^+(H)$. Here is a short proof:
call $h=|G+H|^{-1}$.  ${\rm dim}^+(G+H) = 1+h \sum_{x \in G} {\rm dim}^+ S_{G+H}(x) 
                        +h \sum_{y \in H} {\rm dim}^+ S_{G+H}(y)$.
This is 
$ 1+ h \sum_{x \in G} {\rm dim}^+ (S_G(x) + H) 
     + |G+H|^{-1} \sum_{y \in H} {\rm dim}^+ (G + S_H(y))$.
Induction gets this to
$ 1+ h \sum_{x \in G} {\rm dim}^+ (S_G(x)) + {\rm dim}^+(H)) 
     + h \sum_{y \in H} {\rm dim}^+ (G) + {\rm dim}^+(S_H(y)))$.
This reduces further to ${\rm dim}^+(G) + {\rm dim}^+(H)$.
\index{augmented dimension}
\index{Betre-Salinger formula}

\paragraph{}
There is also the notion of {\bf topological dimension} which is natural if one
looks at the {\bf finite topology of a graph} and which is discussed more
in the next section. The topological dimension is not only a discrete analog of the 
{\bf Lebesgue covering dimension}, but it has the same definition, just used for
finite topological spaces. The classical definition uses the notion of refinement:
a {\bf refinement of a cover} $\mathcal{U}$ is a new cover $\mathcal{V}$ such that
every $V  \in \mathcal{V}$ is contained in some $U \in \mathcal{U}$. A
sub-cover of a cover is an example of a refinement. 
The {\bf topological dimension} now is the smallest $n \geq 0$ such that for
every open cover $\mathcal{U}$ of $G$, there is
a refinement such that every vertex is in no more than $n+1$ sets of the covering.
\index{refinement of a cover}
\index{topology of a graph}
\index{dimension ! topological}
\index{topological dimension}
\index{Lebesgue covering dimension}
\index{dimension ! Lebesgue covering} 

\paragraph{}
The classical topological dimension of the Euclidean space with the topology coming from 
the Euclidean distance metric is $d$ because for every cover, there is a 
sub-cover such that only $d+1$ sets intersect. For $\mathbb{R}$ for example, we
can look at covers by open intervals. We can now constructively get a cover in
which no point is covered by more than 2 intervals. As for an example in the discrete,
the topological dimension of the complete graph $G=K_{n+1}$ is $n$. Any cover has a
refinement in which the sets are the stars $U(v_i)$ where $v_i$ are the vertices 
of $G$. This cover consists of $n+1$ sets.

\paragraph{}
To motivate the next section, let us point out that if we take a finite topological
space coming from a metric, then there are open covers of order $1$ and the
topological dimension is zero. That for a finite topology coming from a metric
space, the topological dimension is always zero can also be seen from the fact that
it produces discrete topology where all sets are {\bf clopen}, open and closed.
In order to have a reasonable topology in which connectivity and dimension works as
expected, we need topologies which are non-Hausdorff. 
\index{topological dimension}
\index{subcover}

\section{Topology}

\paragraph{}
The {\bf topology of a graph} $G$ is a classical {\bf finite topology}
on the simplicial complex $G$ of $G$. The collection of 
{\bf stars} $U(x) = \{ y \in G, x \subset y \}$ together with the
empty set $\emptyset$ define the {\bf basis of this topology}. The basis is
already closed under intersections. Every 
open set in the topology can be written as a union of elements in the basis. 
\footnote{A set of sets closed under intersection is also called a $\pi$-system.}
The {\bf topology} $\mathcal{O}$ is the closure of this basis under the 
operation of taking unions and intersections. The topology is {\bf Alexandrov} 
\cite{Alexandroff1937,May2008},
meaning that any set contains a minimal non-empty open subset. It is not Hausdorff in general 
however. An {\bf open cover} is a finite set of open sets whose union is all of $G$.
An open cover of $G$ necessarily must contain all stars $U(v)$ with $v \in V$. 
The stars $U(v)$ play the role of the ``atoms of space". 
\index{Alexandrov}
\index{topology} 
\index{cover}
\index{refinement} 
\index{topology of a graph} 
\index{open cover}
\index{Alexandrov topology}
\index{star}
\index{$\pi$-system}

\paragraph{}
The topology of a graph has parallels to the {\bf Zariski topology} in algebraic 
geometry, where the closed sets of an algebraic set are the algebraic subsets. 
The closed sets of the graph topology are the sub-simplicial complexes. A slightly 
rougher topology only uses the closed sets which come from sub-graphs of $G$. 
This has some disadvantages like already that for Euler characteristic $\chi(A \cup B)
\neq \chi(A) + \chi(B) - \chi(A \cap B)$ holds only for sub-simplicial complexes and 
not for graphs as the union of two graphs can generate simplices of larger
dimension. $K_3$ as a union of $K_2$ and a path graph $P_3$ intersecting in 
a zero-dimensional $\overline{K_2}$. But the union of the simplicial complexes of $K_2$
and $P_3$ is the simplicial complex of $C_3$ which is not generated by a graph. 
\index{Zariski topology}

\begin{lemma}
The topological dimension is the same than the maximal dimension.
\end{lemma}
\begin{proof}
We follow a proof given in \cite{Gao2024}: 
the topological dimension is less or equal than the maximal dimension $n$ because
for {\bf any open cover} $\{U_j\}$ , the set of all stars $\{U(v)\}_{v \in V}$ of the vertices of $G$ 
is always a open refinement of $\{U_j\}$ because for all $v\in V$ the star $U(v)$ 
must be in one of the open sets $U_j$. There is no $x \in G$ with more than $n+1$ vertices. 
The order of the open cover $\{U(v)\}_{v \in V}$ is therefore $n+1$ or less. 
We have shown that every open cover has a subcover of order $n+1$.
(ii) to see that the topological dimension is larger or equal than the maximal dimension $n$
note that the vertex cover $\{U(v)\}_{v \in V}$ of $G$ of order $n+1$ because $(n+1)$ 
sets are needed to cover $K_{n+1}$, the largest complete subgraph of $G$. 
\end{proof}

\paragraph{}
In general, a basis for the set of closed sets are the stars
{\bf stable manifolds} $U(x) = \{ y \in G, x \subset y \}$ of simplices $x$. These are subsets
of the simplicial complex. One can see a star also as the {\bf unstable manifold} of a simplex $x$.
The {\bf core} $C(x) = \{ y \in G, y \subset x \}$ is a sub simplicial complex which can also be 
seen as the {\bf stable manifold}. 
\footnote{We need to distinguish $x$, which is a single element in the simplicial complex $G$ and 
$W^-(x)$ which is a sub-simplicial complex of $G$.}
Seen in this way, every simplicial complex and in particular, every graph has a {\bf hyperbolic
structure}. 
\index{star}
\index{core}
\index{stable manifold}
\index{unstable manifold}
\index{hyperbolic structure}

\paragraph{}
Not all subsets of the set of simplices are open or closed. For the cyclic graph $C_5$
for example, only 12 percent of all set of simplices are open. And for $K_4$, only
half a percent of all subsets are open. Cohomology can be extended from closed sets to open sets
as they are delta sets
\cite{FiniteTopology,DiscreteAlgebraicSets}. We considered in October 2016
the topology on the simplicial complex of a graph, where the closed sets 
come from sub-graphs which are generated by 
their vertex sets. Not all simplicial sub-complexes are closed in this topology so that
this is a rougher topology. 
The topology in which all simplicial subcomplexes define closed sets is slightly finer.
The topology having subgraphs as closed sets does not work well with counting
because the union of two subgraphs might generate a larger complex than the 
union of the complexes.
\index{Zariski topology}
\index{Cohomology}

\paragraph{}
With a topology comes a {\bf Borel $\sigma$-algebra}
obtained by closing $\mathcal{O}$ under taking complements, unions and intersections. 
A functional like $\mu_k(A) = f_k(A)/f_k(G)$ is now an example of a {\bf probability measure}
on this $\sigma$-algebra. {\bf Measurable sets} are more general objects than graphs
or simplicial complexes. 
Open sets were relevant already in in {\bf connection calculus} because
$g(x,y) = \omega(x) \omega(y) \chi(U(x) \cap U(y))$ were the {\bf Green function entries}
of the matrix $g=L^{-1}$, where $L(x,y) = 1$ if $x \cap y$ is non-empty and $L(x,y)=0$ 
else. We wrote $L(x,y)$ also as the Euler characteristic of a closed set 
$\chi(W^-(x) \cap W^-(y))$. See \cite{KnillEnergy2020,EnergizedSimplicialComplexes2,
  EnergizedSimplicialComplexes3,Sphereformula}.
Any function $h: G \to R$ on a ring extends to a $R$-valued measure on
the $\sigma$ algebra by defining $h(A) = \sum_{x \in A} h(x)$. An {\bf energized 
simplicial complex} defines then a measure on this algebra. 
\index{energized complex}
\index{Borel $\sigma$ algebra}
\index{Green function}
\index{connection calculus}
\index{probability space}

\paragraph{}
The {\bf stable manifold} $W^-(x)=\{ y \subset x \}$ of a simplex $x \in G$ is 
closed because it is a sub-complex the full complex $G$. 
Also its {\bf boundary} $S^-(x) = \delta W^-(x)$, as the intersection of $W^-(x)$ and 
the complement of the open set $U(x)$, is closed. 
The boundary $S^-(x)$ of $W^-(x)$ is always a sphere of dimension $1$ less
than $x$. One can check this by induction with respect to dimension.
The {\bf unstable manifold} $W^+(x)$ of $x \in G$ is open by definition, 
because it is a star and stars have been defined as basis elements generating the topology. 
The open set $U(x)=W^+(x)$ has a boundary a closed set
$S^+(x)$ that is not always a sphere. We write often just $S(x)$ for this closed set 
in $G$ if we work in a simplicial complex. 
But in graph theory, if $v$ is a vertex in a graph, we write $S(v)$ for the 
{\bf unit sphere in the graph} which is the sub-graph induced by the set of vertices attached to $v$. 
The Whitney complex of the graph unit sphere $S(v)$ is isomorphic to the simplicial complex
unit sphere $S(x)$ in the case $x=\{v\}$ has zero dimension.
\index{Stable manifold}
\index{Boundary}
\index{unit sphere}
\index{stable sphere}

\paragraph{}
If the maximal dimension is positive, the topology of a graph is {\bf non-Hausdorff}: 
we can not separate two non-maximal simplices $x,y$ for example. 
\footnote{A simplex is called non-maximal,
if it is contained in a larger simplex}. The topology of a graph has properties
we like: for example, the topology is {\bf connected} if and only if the 
graph is connected. Note that the topology on $V$ coming from the 
{\bf geodesic metric} on a graph would be completely disconnected, 
because every point $\{v\}$ is both open and closed.  An open set $U$ is 
declared to be {\bf contractible} if the graph generated 
by the union of all vertices which are contained in simplices $x_j \in U$ is contractible. 
Alternatively, we could state that a complex $G$ is contractible if its graph
$G=(\mathcal{V},\mathcal{E})$ with $\mathcal{V}=G$ and $\mathcal{E} = \{ (x,y), x \subset y$
or  $y \subset x \}$ is contractible. 
\index{Hausdorff}

\paragraph{}
A finite set $\{ H_j \}$ of sub-graphs $H_j$ of $G$
is a {\bf graph cover} of a finite simple graph $G$,
if every simplex $x$ in $G$ is a sub-graph of at least one of
the graphs $H_j$ and every of the graphs $H_j$ is contractible. The notion of a graph cover is 
graph theoretical and does not require the structure of simplicial complex. 
But the notion is a bit more tricky. It depends on the notion of contractibility for 
open sets. We could define an open set to be contractible, if its closure is contractible. 
We could also define an open set to be contractible if the subgraph generated by its vertices
is contractible. In both cases, 
every open cover $U_j$ of $G$ with contractible $U_j$ defines a graph cover, 
where $H_j$ is the graph generating the smallest closed set containing $U_j$. 
On the other hand, a graph cover defines the open cover $U_j = \bigcup_{v \in H_j} B(v)$. 

\paragraph{}
From the point of view of covering graphs with contractible sets, 
it does not matter whether we deal with {\bf graph covers} or with {\bf open covers}. 
The minimal number of contractible sub-graphs 
that cover $G$ is the same than the number of contractible open sets which cover $G$. There
is no proof required if we define an open set to be contractible if its closure is 
contractible. Graph covers are more intuitive while open covers are closer to what 
we deal with in the continuum. 

\paragraph{}
In the continuum, one sometimes covers a space also with closed sets. This
happened in particular in Lusternik-Schnirelmann theory. 
\cite{Fox} changed it to open sets. Most of the time, 
like \cite{Spanier} (page 279), one sees closed sets in the definition of 
Lusternik-Schnirelmann category.

\index{Lusternik-Schnirelmann category}
\index{graph covers}
\index{open covers}

\paragraph{}
Having a topology associated to a graph, we can define a map between 
two graphs $f: V(G) \to V(H)$ 
to be {\bf continuous}, if the inverse of an open set is open. 
This is equivalent to that the inverse of a closed set is closed. 
Continuous maps are by no means the same than graph homomorphism
and also not the same than simplicial maps as both require vertices to 
map vertices into vertices. A constant map $f(G)=c$ for example is not a 
graph homomorphism (which by definition must map edges to edges),
but it is continuous map from $G \to G$. 

\paragraph{}
A continuous map on graphs also does not define a {\bf simplicial map} of their 
Whitney complexes because for a simplicial map, simplices must be mapped into simplices. 
Simplicial maps map vertices to vertices so that not all continuous maps on the 
topology of $G$ are simplicial maps. 
A simplicial map maps sub-simplicial complexes into simplicial complexes and so, 
open sets into open sets. One calls this an ``open map". 
A simplicial isomorphism is also called a homemorphism. 
If we have a simplicial map, it must preserves order. This means the inverse of a 
closed set must be closed. In other words, it must be continuous. 
\index{continuous map}
\index{simplicial map}

\paragraph{}
The {\bf Barycentric refinement} $G_1$ of a graph $G$ is a new graph for which the
vertices $V_1$ are the elements in the simplicial complex associated to $G$.
The edges are given by connecting two $(x,y)$ if either $x$ 
is contained in $y$ or $y$ is contained in $x$. 
The Barycentric refinement of $G_1$ is denoted $G_2$ etc. The Barycentric refinement is defined
for much more general spaces like cell complexes or delta sets and produces a graph where
the elements are the nodes and two are connected if one is contained in the other. 
(For delta sets, the Barycentric refinement $G_1$ can have other topological properties than $G$. 
Take an open set $U=\{x\}$ in a simplicial complex for example, where $x$ is a facet, a maximal simplex. 
Now $U$ is a delta set by itself. Its Barycentric refinement would be the 1-point graph. 
The Barycentric refinement construction however shows that graphs not much of a
loss of generality. If one does not mind taking a refinement, one can stay within graph theory. 
\index{Barycentric refinement}

\paragraph{}
In the following, we look at a notion of homeomorphic that is motivated by traditional 
notion of homeomorphism of graphs. Two one-dimensional simplicial complexes are called 
{\bf homeomorphic}, if they have a common Barycentric refinement. This notion is used in 
topological graph theory but does not take into account 2 or higher dimensional structures
in the graph. A more general notion of homeomorphism for graphs is by allowing single edge 
refinements (which are by nature local refinements). Every Barycentric refinement in one dimensions
is a composition of local edge refinements. The graph $C_5$ is obtained from $C_4$ by one 
edge refinement. The graph $C_{8}$ is the Barycentric refinement of $C_4$. 

\paragraph{}
In {\bf topological graph theory}, where graphs are embedded in some $2$-dimensional 
surface, the notion of homeomorphism for graphs is enough because the connected components 
on the surface can then also be considered homeomorphic as all faces (2-dimensional cells)
are silently assumed to be balls. Topological graph theory goes beyond combinatorics
because Euclidean spaces and so the concept of infinity is involved. 
We are looking for a notion of homeomorphism
of graphs or finite abstract simplicial complexes which does not involve geometric 
realizations or infinity and which also can constructively be checked in a 
reasonable amount of time. The notion should be effective. 
\index{homeomorphism}
\index{topological graph theory}

\paragraph{}
Here is a variant of a proposal from \cite{KnillTopology}.
Two graphs $G,H$ are declared to be {\bf pre-homeomorphic} if there exists a continuous 
map from some Barycentric refinement $G_n$ to $H$ and a continuous map from some
Barycentric refinement $H_n$ to $G$. By definition, a graph $G$ is pre-homeomorphic to its
$n$'th Barycentric refinement $G_n$.
When considering one-dimensional cases only (triangle free graphs) or graphs with a 1-dimensional 
skeleton complex imposed, then the notion of pre-homeomorphic is what homeomorphic means in 
the topological graph theory literature like for example in \cite{TuckerGross}. 

\paragraph{}
In \cite{KnillTopology} we made the assumption a bit stronger. 
Because we can not prove yet that a graph pre-homeomorphic to a manifold must be a manifold,
we declared $f$ to be a {\bf homeomorphism} if it is a pre-homeomorphism
and additionally, for every atom $U(y)$ of a maximal simplex $y \in H_n$, the open set
$f^{-1}(U(y))$ is an open ball and for every atom $U(x)$ of a maximal simplex $x \in G_n$
the open set $f(U(x))$ is an open ball. With this definition, homeomorphic graphs in which 
one is a $d$-manifold forces also the other to be a manifold. 
\index{homeomorphic graphs} 

\paragraph{}
The {\bf \v{C}ech graph} of the basis $\mathcal{B}$ of the topology $\mathcal{O}$ of a graph $G$
is a new graph which in many cases has the same topological features than the graph itself. 
The {\bf \v{C}ech graph} of the cover defined by all $U(v)=S^+(\{ v \})$ with $v \in V$ is 
the graph $G$ itself. So, $G$ itself can be seen \v{C}ech graph of an open cover
of $G$. The {\bf order} of an open cover is the smallest $m$ such that each point in 
the space belongs to at most $m+1$ open sets in the cover.
The {\bf topological dimension} is the smallest $k \geq 0$ such that for
every open cover $\mathcal{U}$ of $G$, there is a refinement of order $k$.
We have seen that the topological dimension is the maximal dimension. 
\index{\v{C}ech graph}

\section{Category}

\paragraph{}
The {\bf Lusternik-Schnirelmann category} ${\rm cat}(G)$ of a 
graph $G$ is defined as the minimal cardinality $k$ of a 
graph cover $\{ U_j\}_{j=1}^k$ of $G$ using contractible sub 
graphs $U_j$ of $G$. Note that the sub-graphs $U_j$ are not required
to generate themselves within $G$. There are many contractible sub-graphs of a given 
graph that do not generate themselves.
For example, every {\bf spanning tree} of $G$ covers all the vertices of $G$ and is 
contractible. It does generate the entire graph $G$ but a spanning tree 
does not cover the entire graph in general as there are edges which might be missing. 
We want the union of the $U_j$ to cover not only all the vertices but also all 
the edges. 
\index{Lusternik-Schnirelmann category}

\paragraph{}
The Lusternik-Schnirelmann category is of interest in graph theory also
because it is related in spirit to a class of functionals which are classical and well studied 
in graph theory. For example, the {\bf edge arboricity} tells how many forests are 
needed to cover a graph.  By the {\bf Nash-Williams theorem}
\cite{NashWilliams1964}, this arboricity of a graph $(V,E)$ 
is the maximum of $|E_W|/|W|$, where $W$ ranges over all subsets of $V$
and $E_W$ is the number of edges generated by the induced graph of $W$. 
The {\bf vertex arboricity} (=point arboricity) is the maximal number of forests 
partitioning $V$ such that each forest generates itself. The vertex arboricity 
is of interest because the {\bf chromatic number} is sandwiched between 
the vertex arboricity and twice the vertex arboricity. See 
\cite{ThreeTreeTheorem,ArboricityManifolds}.
\index{Vertex arboricity}
\index{Edge arboricity} 
\index{Chromatic number}
\index{Nash-Williams}

\paragraph{}
Instead of graph covers, we could also use (like Fox classically did in 
\cite{Fox}) use {\bf open covers} of the topology and the category would 
be the same: for an open set $U$, we can look at the graph generated by the 
vertices appearing in $U$ and define $U$ to be contractible, if that graph 
is contractible.  With this definition, we also could have used {\bf closed 
coverings} $U$, the reason being that the closure of a contractible open set 
must by definition be contractible. Unlike for closed sets, where the notion 
of contractible is intrinsic and does not depend on where it is embedded 
into, the notion of contractible for open 
sets depends on the embedding as we look at the closure of the open set. 
We chose here the more graph theoretical definition and use covers by 
subgraphs. 
\index{Open cover}
\index{Closed cover}

\paragraph{}
The topological closure of an open set is not necessarily contractible
as the example of a punctured sphere = open ball shows. An open ball should
with this definition not be considered contractible. 
Contractibility for open sets is a bit tricky is because the cohomology of 
an open set is in general much different than the cohomology of a closed set. 
An open ball has only a non-trivial maximal cohomology. 
Nevertheless, if we wanted to define Lusternik-Schnirellmann category 
using open sets, we would have to define an open set
to be contractible, if its closure is contractible.
\index{Lusternik-Schnirelmann category}
\index{contractible graph cover}
\index{topological closure}

\paragraph{}
By induction, a contractible graph $G$ has Euler characteristic $1$.
The Lusternik-Schnirelmann category of a contractible graph by definition is equal to 
$1$. A $d$-sphere $G$ of dimension $d \geq 0$ almost by 
definition has Euler characteristic $1+(-1)^d$ because removing an open $d$-simplex
of Euler characteristic $(-1)^d$ produces a closed ball of Euler characteristic $1$. 
A $d$-sphere almost by definition also has category $2$ because 
the unit ball $B(v)$ together with the graph $G \setminus v$ cover $G$. Now, the  
unit ball $B(v)$ and $G \setminus v$ are contractible and cover the graph. 

\paragraph{}
Also, almost by definition, the Lusternik-Schnirelmann category of a non-connected graph 
$G$ is the sum of the categories of its connected components. Because 
a discrete manifold is never contractible - we can not even remove a single point - 
the category of a manifold is always $\geq 2$, with the single exception of the 
{\bf $1$-point graph} $K_1$ which is the only connected $0$-dimensional 
manifold. The $1$-point graph $K_1$ is the {\bf terminal object} in the category of graphs. 
\index{terminal object}

\paragraph{}
Lets look at some examples: \\
{\bf 1)} ${\rm cat}(K_n)=1$ for all $n \geq 1$. \\
{\bf 2)} ${\rm cat}(C_n)=2$ for all $n  \geq 4$.  \\
{\bf 3)} ${\rm cat}(T)=1$ if $T$ is a tree as a tree is contractible. \\
{\bf 4)} ${\rm cat}(F)=b_0(F)$ if $F$ is a {\bf forest} a disjoint union of trees. 
   The cover is given by the maximal subtrees. \\
{\bf 5)} The graph complement $\overline{G}$ can have a very different Lusternik-Schnirelmann 
   category than $G$.  For $K_n$, the graph $\overline{K_n}$ has no edges an 
   ${\rm cat}(K_n)=1, {\rm cat}(\overline{K_n})=n$. \\
{\bf 6)} For any d-sphere, ${\rm cat}(G)=2$ as $B(x), G \setminus x$ are both balls and
   so contractible. 

\section{Operations} 

\paragraph{}
The {\bf Stanley-Reisner product} $G \times_1 H$ of two graphs can be defined as the Barycentric 
refinement of the Cartesian product $G \times H$ of its simplicial complexes. It is best defined
in terms of graphs. Let $V=G \times H$ be the vertex set of a graph and $E$ the set of pairs (a,b)
for which either $a \subset b$ or $b \subset a$. The Whitney complex of this graph is now declared
to be the {\bf Stanley-Reisner product} or topological product. This product is not 
associative because $G \times 1=G_1$ is the Barycentric refinement of $G$. 
\index{Stanley-Reisner}
\index{topological product} 

\paragraph{}
The {\bf Stanley-Reisner picture} writes the product in terms of a multiplication of polynomials. 
We have an associative product within polynomials. But this does not mean associativity 
for the product. The reason is that we need to returning from a ring element to a simplicial complex. 
If $G$ is a complex and $p$ the Stanley-Reisner polynomial then the complex of $p$ is 
the Barycentric refinement $G_1$. 

\paragraph{}
Since the Cartesian product of two simplicial complexes (as sets) is not a 
simplicial complex in general, we must use a product in 
a larger category of cell complexes or even the larger {\bf topos of delta sets} 
The later is a natural frame work to do finite geometry \cite{EilenbergZilber}. 
The Barycentric refinement always could bring us back to graphs. 
The upshot is that a Cartesian product for graphs which works in higher dimension and
satisfies the dimension and Kuenneth requirements can be done either with the Stanley-Reisner
product (no associativity), abstractly by extending the frame work or by using delta sets. 
Both the Cartesian product of graphs or the categorial product (small product) of graphs are 
not suited for higher dimensional considerations. A good product in graph theory is the
Shannon product, but it does not preserve manifolds. 

\paragraph{}
If we want to use graphs and stay in the category of graphs and have associativity, we need to 
pick one of the few products available: small product, Categorical product, large product (Sabidussi)
and strong product (Shannon). None of these products really fits all bills like that the product of
manifolds is manifolds. The Shannon product has many nice properties, the Shannon product 
of two graphs is homotopic to the topological Stanley-Reisner product.
The Shannon product together with disjoint union as addition defines a semi-ring that 
extends to a Shannon ring of graphs. But the Shannon product of 
two manifolds is only homotopic to a manifold and not a manifold by itself in general. 
\index{topos}
\index{delta set}
\index{Stanley-Reisner}
\index{Topological product}
\index{Shannon ring}

\paragraph{}
Classically, the Cartesian product of two (Euclidean) $k$ manifolds without boundaries 
is known to have category at least $(k+1)$. We will show this later. 
A {\bf k-torus} for example has category $k+1$. 
Category is {\bf not a homotopy invariant}: the {\bf dunce hat} $G$ for
example is homotopic to $1=K_1$ but $G$ has category $2$ as the dunce hat is not 
contractible. We could define the homotopy invariant
$\overline{cat}(G)$ as the smallest category of a space homotopic to $G$.
\index{k-torus}
\index{homotopy invariant}
\index{cup product}

\footnote{
If we would define $\overline{cri}(G)$ is the minimum of 
all ${\rm cri}(H)$ with $H$ homotopic to $G$, the main theorem of 
Lusternik-Schnirelmann theory would then mean
${\rm cup}(G)+1 \leq \overline{{\rm cat}}(G) \leq \overline{{\rm cri}}(G)$
relating three homotopy invariants ${\rm cup} = \overline{{\rm cup}}$, $\overline{{\rm cat}}$
and $\overline{{\rm cri}}$.  }
\index{category ! homotopy invariant}

\paragraph{}
In the following, we mean with the product the Cartesian product. 
The {\bf join} (or {\bf Zykov join}) 
of two graphs $G,H$ is the disjoint union with additional connections 
between any vertex of $V(G)$ and $V(H)$. For simplicial complexes $G,H$
the join is the disjoint union $G \cup H$ together will all elements
$\{ x \cup y, x \in G, y \in H \}$.  
\index{join ! graph}
\index{join ! simplicial complex}
\index{Zykov join}

\paragraph{}
We denote with $G * H$ the {\bf Shannon product} of $G$ and $H$. It is 
the graph for which the vertex set is the Cartesian product of the 
vertex sets of $G$ and $H$ and where $(x,y)$ is connected with $(a,b)$
if one of the projections is an edge in $G$ or $H$. 
\index{Shannon product} 

\begin{lemma}
a) ${\rm cat}(G_1)={\rm cat}(G)$ if $G_1$ is the Barycentric refinement. \\
b) For the disjoint union $+$, ${\rm cat}(G+H) = {\rm cat}(G)+{\rm cat}(H)$. \\
c) The join $G \oplus H = \overline{\overline{G} + \overline{H}}$ gives
   ${\rm cat}(G \oplus H) = {\rm min}({\rm cat}(G), {\rm cat}(H))$. \\
d) ${\rm cat}(G \times H) \leq {\rm cat}(G) \cdot {\rm cat}(H)$. 
e) ${\rm cat}(G*H)        \leq {\rm cat}(G) \cdot {\rm cat}(H)$. 
\end{lemma}
\index{minimal number of critical points} 

\begin{proof} 
 a) The Barycentric refinement of a contractible graph is contractible.
    If $\{ U_j \}$ is a cover of $G$ by contractible graphs then the 
    Barycentric refinements $\{ V_j = (U_j)_1 \}$ form a cover of $G_1$
    by contractible graphs. \\
 b) If $\{ U_i \}$ is a cover of $G$ and $\{ V_j\}$ is a cover of $H$, 
    then $\{ U_i \} \cup \{ V_j\}$ is an open cover of $G \cup H$. We can not take less. \\
 c) If $U$ is a contractible graph and $H$ is an arbitrary graph, then 
    $U \oplus H$ is contractible. Assume that $G$ is the graph with minimal
    category $k$. This means we have a cover $U_1, U_2, \dots, U_k$ of $G$. 
    The graphs $U_j \oplus H$ are contractible. The graph $G \oplus H$ therefore
    can be covered with $k$ contractible graphs. \\
 d) If $\{ U_i \}$ is a cover of $G$ and $\{ V_j \}$ is a cover of $H$ then 
    $\{ U_i \times_1 V_j \}$ is a cover of $G \times H$.  \\
 e) If $\{ U_i \}$ is a cover of $G$ and $\{ V_j \}$ is a cover of $H$ then 
    $\{ U_i * V_j \}$ is a cover of $G \times H$ and
    $U_i * V_j$ are contractible.
\end{proof} 

\paragraph{}
The example of the torus $\mathbb{T}^2 = \mathbb{T}^1 \times \mathbb{T}^1$ shows
that part d) can not be an equality in general because ${\rm cat}(\mathbb{T}^2)=3$. 
However, if ${\rm cat}(G)=1$, then ${\rm cat}(G \times H) = {\rm cat}(H)$. 

\paragraph{}
Let us add a remark to e). Since 
${\rm cup}(A*B) \geq {\rm cup}(A) {\rm cup}(B)$ 
we have ${\rm cat}(A*B) \geq {\rm cup}(A*B)+1 \geq {\rm cup}(A) {\rm cup}(B)$
For example, the category of a Shannon product of a non-simply connected orientable
manifold is $\geq k+1$ because ${\rm cup}(G)+1 \leq {\rm cat}(G)$ and the
fact that every discrete manifold that is orientable and non-simply connected has
a degree $1$ differential form. 

\section{Functions}  

\paragraph{}
When looking at functions on a graph, we want to impose some non-degeneracy condition
in general. A most natural one is that if $(x,y)$ is an edge then $f(x) \neq f(y)$.
One calls this {\bf locally injective} or a {\bf coloring}.
This local continuity condition was needed in Poincar\'e-Hopf. In the continuum, one imposes more
regularity on maps by assuming that the functions to be Morse, meaning that at
all critical points, the Hessian is invertible. In the continuum, one also works
in a smooth setting, which means that at a critical point $x$, 
the stable spheres $S_r^-(x) = S_r(x) \cap \{ f \leq f(x) \}$ are spheres. 
This is what one can adopt also in the discrete. At regular points, the 
stable sphere $S_r^-(x)$ is a ball for small enough $r$. We adapt this to the discrete.
\index{coloring}
\index{locally injective} 

\paragraph{}
For a general finite simple graph $G=(V,E)$ and a locally injective
function $f$, a vertex $x \in V$ is called a {\bf regular point}
if $S^-_f(x)$ is contractible. If $S^-_f(x)$ is not contractible,
it is called a {\bf critical point}. 
The minimal number of critical points of a general locally injective function
is denoted ${\rm cri}(G)$.
\index{regular point}
\index{critical point}

\paragraph{}
As defined in the introduction already, 
a function $f:V \to \mathbb{R}$ is a {\bf Morse function} 
if every $S^-(x)$ for $x \in V$ is either a $k$-sphere for some $-1 \leq k \leq n-1$ 
or contractible. If $S^{-}(v)$ is a contractible, the 
point is called a {\bf regular point}, otherwise it is considered to be a
{\bf critical point} of $f$. A function could be called {\bf strongly Morse} 
if both $f$ and $-f$ are Morse. The later condition is hard to achieve if we have 
not a manifold. It does not required it however: take a manifold $M$ and make a homotopy 
extension by adding one vertex $w$ to a single vertex of $M$. If $g$ was a strongly Morse
function on $M$, then extend it to the additional point by setting 
$g(w)= {\rm max}_{v \in V(M)} g(v) + 1$. 

\paragraph{}
An example of a strongly Morse function
is the {\bf dimension function} on the Barycentric refinement $G_1$ of a graph $G$. 
If $v$ is a critical point of a Morse function, then $S^-f(x)$ is a $(k-1)$-sphere
for some $k$ and $k$ is called the {\bf Morse index} of $x$. 
Because of the {\bf Euler-Gem formula} $\chi(A)=1+(-1)^k$ for any 
$k$-sphere $A$, the Poinar\'e-Hopf index of a point $v$ of Morse index $k$ is
$i_g(x) = (-1)^{k}$. If $c_{k}$ is the number of critical 
points of Morse index $k$ in $G$, then the {\bf Poincar\'e-Hopf} theorem
reads as $\sum_k (-1)^k c_k  = \chi(G)$, where
$\chi(G) = \sum_k (-1)^k f_k(G)$, and $f_k(G)$ counts the 
number of $k$-dimensional simplices. 
\index{Morse function}
\index{Euler Gem}
\index{Critical points}

\paragraph{}
For a Morse function $f$ on a discrete manifold, the general critical points
agree with critical points of the Morse function because for a Morse function,
all $S^-(x)$ are either a ball and so contractible or spheres and
so non-contractible. The minimal number of critical points which some Morse function can
achieve on a discrete manifold is denoted by $c(G)$.
We have ${\rm cri}(G) \leq c(G)$. Proof: for a Morse function $f$ with
$c_f$ critical points, ${\rm cri}(G) \leq c_f(G)$.
Now minimize the right hand side over all Morse functions $f$ to get
${\rm cri}(G) \leq c(G)$. 

\paragraph{}
The inequality can be strict like for any 2-torus graph
$\mathbb{T}^2$, where ${\rm cri}(G)=3$ and $c(G)=4$.
Indeed, for a Morse function, the indices are $\pm 1$ and by Poincar\'e-Hopf
they add up to $0$ so that there must be an even number of critical points for 
a Morse function on a $2$-torus. There can not be $2$ because that would 
mean that it is a sphere. 

\index{Exterior derivative} 

\paragraph{}
The vertex sets $G$ of all complete sub-graphs of $(V,E)$ is a
{\bf finite abstract simplicial complex}, a finite set of non-empty 
sets closed under the operation of taking non-empty subsets. 
The finite structure was introduced in 1907 by Dehn and Heegaard 
\cite{DehnHeegaard,BurdeZieschang}.
We can equip each of these sets $x$ with an {\bf orientation}, a fixed order
of its vertices. An example is to enumerate the vertices of $G$, then use
these labels on each simplex. The simplicial complex of the kite graph 
$K_{1,2,1}=(V,E)=(\{1,2,3,4,5\},\{ (12),(13),(23),(24),(34) \})$ for example is
$G=\{ \{1\},\{2\},\{3\},\{4\},\{5\},\{1,2\},\{1,3\},\{2,3\},\{2,4\},\{3,4\},
\{1,2,3\},\{2,3,4\} \}$. Writing this down in lexicographic order already gives
an example of an orientation. 

\paragraph{}
A {\bf $k$-form} is a function on $G_k$,
the set of sets in $G$ with $k+1$ elements. \footnote{The function could be 
symmetrized by defining $f( \pi(x)) = {\rm sign}(\pi) f(x)$
if $\pi$ is a permutation.  It would then be anti-symmetric but there is no need to do that 
and just consider $f$ as a function on the ordered element.}
A differential form $f$ is nothing else than a scalar function on $G = \cup_{k=0} G_k$, 
once the elements of $G$ are equipped with orientations. The restriction of $f$ to $G_k$
is the vector space of {\bf k-forms}. The dimension of this vector space is $f_k$. 

\paragraph{}
The total set of forms is a $n=\sum_{k=0}^d f_k$ -dimensional vector space $\Omega$, 
where $n$ is the number of simplices of $G$ and $f_k=f_k(G)$ the number of $k$-dimensional 
simplices in $G$. The set of $k$-forms $\Omega_k$ is a $f_k$-dimensional real vector space. 
Define {\bf the exterior derivative} 
$$ df(x) = \sum_{y \subset x, {\rm dim}(y)={\rm dim}(x)-1} {\rm sign}(x,y) f(y) \; , $$
where ${\rm sign}(x,y)$ is $1$ if the orientation of $y$ matches the orientation
of $x$ on $y$ and $(-1)$ else. \footnote{One usually writes this 
as $df(x) = \sum_{j=1^k} (-1)^{k-1} f(x_1,\dots,\hat{x}_j,\dots,x_k)$ which is the same
because if $y$ is the simplex with the $j$'th vertex removed, then $y$ has the
same orientation than $x$ if and only if $(j-1)$ is even.}
For example, ${\rm sign}( \{1,2,3\},\{1,3\} )=-1$ and ${\rm sign}( \{1,2,3\},\{2,3\} )=1$. 
\index{exterior derivative}

\paragraph{}
If $G$ has $n$ elements, the exterior derivative map $d$ is represented in an explicit way
as a lower triangular $(n \times n)$-matrix satisfying $d^2=0$. 
If $d^*$ is the {\bf transpose} of $d$, the matrix $D=d+d^*$ is 
called the {\bf Dirac operator} of $G$ \cite{KnillILAS} and 
$L=D^2=d d^* + d^* d$ is the {\bf Hodge operator} of $G$. Note that the
matrix entries of both matrices $D$ and $L$ depend on the orientation choice of the 
simplices but that any different choices of the orientation just corresponds to an
orthogonal change of coordinates in the finite dimensional Hilbert space $\mathbb{R}^n$
or $\mathbb{C}^n$ on which $D$ and $L$ operate. The reason for the name ``Dirac operator" 
is because $D$ is a square root of $L$ which shares properties of the Dirac operator in 
the continuum. The operator $D=d+d^*$ also appears in the continuum \cite{Cycon}.

\paragraph{}
The map $d$ maps $\Omega_k$ to $\Omega_{k+1}$ and
$d^*$ maps $\Omega_{k+1}$ to $\Omega_{k}$.
Because $L$ leaves the $f_k$-dimensional space $\Omega_k$ invariant, $L$ decays 
into blocks $L_k= d_k^* d_k + d_k d_k^*$, for $k \geq 0$. Each $L_k$ is a $f_k \times f_k$
matrix. Changing the orientation of a simplex is a unitary transformation.
The non-negative numbers $b_k(G)={\rm dim} {\rm ker}(L_k)$ are called the {\bf Betti numbers}. 
They do not depend on the choice of orientation. 
The kernel ${\rm ker}(L_k)=H^k(G)$ are {\bf k-harmonic forms} and represent 
cohomology classes of $G$. It is isomorphic to the traditionally defined as
the space ${\rm ker}(d_k)/{\rm im}(d_{k-1})$. This is the same because of the {\bf Hodge decomposition}
${\rm im}(d_k) \oplus {\rm im}(d_{k+1}^*) \oplus {\rm ker}(L_k)$ which follows directly from 
the {\bf rank-nullity theorem} applied to the block side diagonal matrix $L_k$ with $d_k,d_{k+1}$ blocks
so that ${\rm im}(d_k) \oplus {\rm im}(d_{k+1}^*)$ is ${\rm im}(L_k)$. 
\index{Hodge operator}
\index{Dirac operator}
\index{Betti numbers}
\index{Cohomology classes} 
\index{Harmonic forms}
\index{Hodge decomposition}

\paragraph{}
The zero'th block $K=L_0$ of the Hodge Laplacian $L$ is 
known as the {\bf Kirchhoff matrix}. It can be written as $B-A$, where $B$ is the diagonal 
{\bf vertex degree matrix} with $B(v,v)$ being the {\bf vertex degree} of the vertex $v$ and 
$A$ is the {\bf adjacency matrix} of $G$ which is a $0-1$ matrix with $A(v,w)=1$ if and only 
if $(v,w) \in E$. From this representation $K=B-A$, one can see that $K$ does not depend 
on the choice of the orientation used on vertices. The definition goes over to 
quivers, graphs for which multiple connections and self-loops are allowed 
\cite{Eigenvaluebounds}. 

\paragraph{}
A small example, where some matrix entries of $L$ 
can depend on the orientation is the kite graph $G=K_{1,2,1}$ obtained 
by taking away a vertex from $K_4$. The off diagonal entries of $L_2$ can depend on the orientation.
But in general, the change of orientation of the simplices only induces an orthogonal change of 
basis in the linear space $\Omega$. It is like choosing a coordinate system in traditional 
Euclidean geometry. 
\index{Kirchhoff matrix}
\index{Hodge Laplacian}
\index{kite graph}

\paragraph{}
The {\bf Shannon product} of two graphs $G=(V,E),H=(W,F)$ has as vertex set the 
{\bf Cartesian product} $V \times W$ and as edge set 
$\{ ((a,b),(c,d)), (a=c)$ or $(a,c) \in E)$ and 
($b=d$ or $(b,d) \in F) \}$. Two pairs in the Cartesian product are connected 
if the projection on each of the factors is either a vertex or edge and 
at least one is an edge. A $k$-form as a function on the set of signed 
$k$-simplices. Given a $p$-form $f$ on $G$ and a $q$-form $g$ on $H$, define 
a $(p+q)$-form on the Shannon product $G*H$ by first considering 
$f*g(x,y) = f(x)g(y)$ which is a $(k+l+1)$-form, then take the divergence.
When restricted to cohomology classes, we get the K\"unneth formula \cite{Kuenneth}
(see \cite{KnillKuenneth} for graphs). 
It is a function on $(k+l)$-simplices on $G*H$ and so is a $(k+l)$ form 
on $G \times H$. This gives the product $f \otimes g = d^*(f*g)$. 
\index{Shannon product}
\index{Kuenneth formula}

\paragraph{}
In order to have an associative product on simplicial complexes, one has to 
leave the category of simplicial complexes. There is no Cartesian product of 
simplicial complexes that is associative, has the natural dimension properties 
and preserves manifolds. One can go beyond simplicial complexes by using 
cell complexes like discrete $CW$ complexes. 
The most convenient one is the even larger class of  {\bf delta sets}. 
Delta sets form a {\bf functor category} and are {\bf presheaves}. It is usually given 
by a finite set of sets $G_i$ and face maps $d_i^n: G_{n+1} \to G_n$ satisfying 
$d_i^n d_j^{n+1} = d_{j-1}^n d_i^{n+1}$ leading to the exterior derivative 
which maps $f \in \Omega_n$ into 
$df(x) = \sum_{i=0} (-1)^i f(d_i^n x) \in \Omega_{n+1}$. The commutation relation 
in the axiom for delta set assures that $d$ is an exterior derivative. 
\index{Functor category}
\index{presheave}
\index{CW complex} 
\index{cell complex} 
\index{delta set} 

\paragraph{}
It is more convenient to encode a delta set as a finite set $G$ with $n$ elements,
a single $n \times n$ matrix $D=d+d^*$ and a dimension function 
$r:G \to  \mathbb{N}$ which selects the cells $G_n = r^{-1}(n)$ of dimension $n$. 
The information $(G,D,r)$ encodes everything we need to know. 
The category of delta sets is powerful because it is a {\bf topos} and so 
Cartesian closed. The same tensor product construction works for differential forms,
but instead of bringing down the dimension, we just declare elements $(x,y)$ as $(p+q)$-dimensional
if $x$ has dimension $p$ and $y$ has dimension $q$. For example, if $x,y$ are both $0$-dimensional, 
then $(x,y)$ is a $0$-dimensional point even, so the data structure represents it as a 
$0$-dimensional object. 
\index{delta set}
\index{topos}
\index{Cartesian closed category}

\section{Cup product}

\paragraph{}
As Hassler Whitney first realized \cite{Whitney1992}, the
definition of the cup product is a bit puzzling at first in a discrete setting. 
If we take a $k$-simplex $x$ and a $l$-simplex $y$, then $x \oplus y$ is a $k+l+1$-simplex. 
For example, if $x=(a,b)$ and $y=(c,d)$ are $1$-dimensional, 
then $x \oplus y = (a,b,c,d)$ is a tetrahedron and $3$ dimensional. We would like however
to have tensor product definition as in the continuum and get a 2-dimensional object. 
The trick is to take pointed simplices and join them along this point. For example, 
if $(a,b,c)$ is a triangle in a graph, then fix a point $a$. A 1-form cocycle  on the 
triangle is now fixed by giving $f(a,b)$ and $f(b,c)$. Given two $1$-form cocycles, 
define $f \otimes g(a,b,c) = f(a,b) g(a,c)$. The wedge product is now 
$f \wedge g(a,b,c) = f(a,b) g(a,c) - f(a,c) g(a,b)$. This parallels what the cross
product does in $\mathbb{R}^2$. The vector space of $1$-cocycles on $(a,b,c)$ 
is two dimensional and can be described by a vector $[x,y]$ where $x=f(a,b), y=f(a,c)$. 
Given two such 1-cocycles $[x_1,y_1],[x_2,y_2]$, the cross product is 
$x_1 y_2 - x_2 y_1$.  \\

\paragraph{}
If we do the same thing on a 3-simplex (a tetrahedron) $(a,b,c,d)$ then the space of $1$-cocycle
is 3-dimensional. It can naturally be identified as $\mathbb{R}^3$ containing elements
$[x,y,z]$ where $x=f(a,b),y=f(a,c),z=f(a,d)$. The exterior product is now a 2-cocycle which 
is determined by the values on the three triangles containing $a$. The 2-cocycles
also form a 3-dimensional vector space $\mathbb{R}^3$ containing elements
$[x,y,z]$ where $x=f(a,c,d), y=f(a,d,b), z=f(a,b,c)$. The exterior product is now the naturally
the cross product after identifying the vector spaces of 1-forms and 2-forms. 
Lets now look at the general case: 

\paragraph{}
If $x$ is a $(k+l)$-simplex in $G$, it can be written as $x=(x_0,x_L,x_R)$ 
Define $f {\rm cup} g = f \otimes g + (-1)^{|x_l|} g \otimes f$
which defines a $(k+l)$-form on $G$. 
For example $x=(x0,x_1,x_2)$ is $1+1$ simplex if $(x_0,x_1)$ and $(x_0,x_2)$ are 1-simplices.
Then, $f \otimes g(x) = f(x_0,x_1) g(x_0,y_1)$ and 
$f \wedge g(x) = f(x_0,x_1) g(x_0,y_1)-g(x_0,x_1) f(x_0,y_1)$. 
Now $d (f \otimes g) = df \otimes g + (-1)^k f \otimes dg$. 

\paragraph{}
The {\bf exterior product} can be defined cocycles of  any simplicial complex $G$. 
It is associative and super commutative so that it is a graded {\bf super algebra}. 
In the continuum, it is also known as the {\bf Grassmann algebra}.
Since cocycles are mapped by $d$ into cocycles and coboundaries into coboundaries, 
different cohomology classes are mapped into different cohomology classes. 
This produces the {\bf cup product} on the space of harmonic forms 
$H^p \times H^q \to H^{p+q}$. 
\index{exterior product} 
\index{Grassmann} 
\index{cup product}

\paragraph{}
There are different ways to define the exterior product in the discrete. 
In \cite{josellisknill}, we chose to symmetrized version of the product 
to make it orientation independent.  Here are some properties of 
the exterior product. 

\begin{lemma}
a) Leibniz rule: $d (f \wedge g) =df \wedge g + (-1)^p f \wedge dg$.  \\
b) Associativity $(f*g)*h =f *(g*h)$ and see that it is a coboundary. \\
c) Super commutativity $f \wedge g = g \wedge f (-1)^{p q}$. \\
d) If $df=0$ and $dg=0$, then $d(f*g)=0$. \\ 
e) If $f=dh$ and $g=dk$ gives $f*g = dh*dk = d (h*dk) = -d(dh*k)$ \\
f) The one-element is the constant harmonic $0$-form which is 1 everywhere. 
\end{lemma}

\paragraph{}
A simplex is {\bf locally maximal} if it is not strictly contained in an other simplex. 
A graph is {\bf orientable} if one can orient the locally maximal simplices in such a way 
that ${\rm sign}(x,y)=1$ for all $x \subset y$ as well as for $y \subset x$. 
Not all graphs are orientable. The smallest non-orientable graph is the graph 
complement of $C_7$ which is a discrete {\bf Moebius strip} with $7$ vertices, $14$ edges 
and $7$ triangles. Fixing an orientation of a triangle fixes the orientation of 
its neighbors but going around the closed loop changes the orientation. 
An orientable complex, we can define a Hodge dual. 
\index{locally maximal}
\index{orientable}
\index{Moebius strip}

\paragraph{}
In an orientable discrete 2-manifold, the exterior product of two 
$1$-forms in different cohomology classes is a $2$-form, a {\bf volume form}. 
In three dimensions, the exterior product defines {\bf cross product} when
identifying 2-forms and 1-forms. In classical calculus, 
where we associate both 1-forms or 2-forms as {\bf vectors}, we can build
a cross product $v \wedge g$ which is considered a scalar function but in more advanced
set-ups like differential geometry the cross product can be related to the exterior
product obtained from two 1-forms.  For example, if $G=K_3$ is oriented cyclically
and $f$ assigns the values $(1,0,0)$ to the three edges and $g$ the values $(0,1,0)$
to the three edges, then $f \wedge g$ gives the value $1$ on the 2-simplex. 
\index{cross product ! two dimensions}
\index{volume form}
\index{cross product}

\paragraph{}
In an orientable $3$-manifold, we also have
{\bf dual edges} $E'$. They are given by pairs $(x,y)$, where $x$ is a 2-simplex and 
$y$ is a $0$-simplex such that $x,y$ span a 3-simplex. 
We can think of a dual edge as an {\bf altitude} from the vertex $y$
to the center of $x$. Since in a 3-manifold every 2-simplex $y$ bounds two 0-simplices $x,z$,
and if the maximal simplices define the orientation, then there is for every triangle a unique
altitude $(x,y)$ such that $y$ is compatible with the orientation of the tetrahedron $x \cup y$. 
A $2$-form therefore can be visualized as attaching values to these dual vectors. If we chose
a coordinate system $(i,j,k = i \times j)$ in each tetrahedron and postulate the relations of
Hamilton $i^2=j^2=k^2=ijk=-1$, the multiplication can be extended from edges $E$ to the union 
of edges and dual edges $E \cup E'$. This means that in every tetrahedron, we have a 
quaternion algebra. The exterior product makes any orientable 3-manifold into a Lie algebra, 
once we join 1 and 2-forms. 

\paragraph{}
The closest to the continuum is to define a {\bf Hodge dual} to a $2$-form $f$ as a
$1$-form $g$ which has the property that $f \wedge g$ is the volume form. The Hodge dual can be 
defined if one has a volume form. The Hodge dual is well defined if to $y \subset x$ is given 
defines a unique $z$ with $y+z=x$ overlapping in the smallest element. Having this, allows
us to define a cross product of 1-forms in a graph given by an orientation. 
\index{Hodge dual}

\paragraph{}
Let $G$ is a $d$-manifold and $f$ is a locally injective function on $G$. 
If $c$ is a value not taken by $f$, 
we can look at the sub-graph $A=\{ f = c \}$ of $G_1$ 
generated by all simplices $x$ of $G$ on which $f-c$ changes sign.  \cite{KnillSard}
That this graph is a $(d-1)$ manifold is proven by showing that
for every simplex $x$, $\{ f = c \} \cap S(x)$ is by induction assumption 
a $(d-1)$ manifold. Locally injective functions which have the property that 
$f(x) = \sum_{v \in x} f(v)/|x|$ produces locally injective functions on the 
Barycentric refinements can now be used to define ``varieties" $\{f_j=c_j, j=1, \dots ,k\}$ 
in the $k$'th Barycentric refinement. These are all vertices in $G_k$ on which all functions 
$f_j-c_j$ change sign. It is always a $(n-k)$-manifold or empty. 

\paragraph{}
For example, If $G$ is a 3-manifold, then at every vertex we can look at
$\{ f = c \} \cap S(x)$ which is a disjoint union of circular graphs. If $G$ is a $4$-manifold,
then for any two functions $f_1,f_2$ we can look at the 3-sphere $S(x)$ and have
a {\bf link} $\{ f_1=c_1, f_2=c_2 \}$, a one-dimensional sub-manifold of the 3-sphere. 
For a given $n$ manifold, almost all functions $f$ we can look at manifolds $f=c$ with 
$min f<c<max f$. We have so a probability space of sub-manifolds. We can look at expectation
values. We could look for example what the average $k$'th Betti number is or what the 
average Euler characteristic is Or we can look how many of the functions are Morse 
functions. 

\section{Morse Inequality}

\paragraph{}
We get now to the Morse inequality. For the classical theory, see 
\cite{Morse1939,Matsumoto,NicoalescuMorse}. A semiclassical analysis approach 
to the Morse inequalities using Witten deformation \cite{Witten198} can be 
found in \cite{Cycon}. For Forman's discrete Morse theory \cite{forman98,Forman2002}. 
Morse theory competes with Lusternik-Schnirelmann theory.
We will later see a comparison of L-S with a Picasso painting. Morse
theory was promoted in  \cite{BottIndomitable} by:
{\it And the term " critical point " of course brings me to my topic proper of this 
morning: "Morse Theory Indomitable ". I think Morse would have approved the title 
for when I first met him, he preached the gospel of critical point theory first, 
last and forever, to such an extent, that we youngsters would wink at each other 
whenever he got started}.

\paragraph{}
Unlike in the continuum, 
the theorem works for any graph, whether it is a manifold or not. It works for example for 
the Barycentric refinement of a simplicial complex $G$ and the function $f(x)={\rm dim}(x)$,
which is a really special case, where $b_k=c_k$, but still where $G$ is not necessarily
a manifold. The result works also for simplicial complexes. 
In that case, the Morse extensions have to be done by adding new simplices by 
joining a given simplicial complex (a closed set) with a star (an open set). 
An example of a homotopy extension $K_2 \to K_3$ adds a new vertex to the contractible
$K_2$. In the simplicial complex picture, we would add an open cone
$U=\{ \{3\},\{1,3\},\{2,3\},\{1,2,3\} \}$ to $\{ \{1 \},\{2\},\{1,2\} \}$. One can then 
see $U=K_3 \setminus K_2$ as the complement of the closed set $K_2$ in $K_3$. \\

{\bf Examples:}  \\
1) The complete multipartite graph $K_{(2,2,2)}$ is the octahedron graph, a 2-sphere. 
   Its Betti vector is $(1,0,1)$ like for all 2-spheres. There are Morse functions on the 2-sphere
   with exactly two critical points. One with index $0$ (the minimum) and one with index $2$,
   the maximum.  \\
2) The complete multipartite graph $K_{(2,2,2,2)}$ is the smallest 3-dimensional sphere. 
    Its Betti vector is $(1,0,0,1)$. Also here, there is always a Morse function which has exactly 
    two critical points. This works in arbitrary dimensions and is related to a theorem of Reeb
    in the continuum \cite{Reeb1952}.  \\
3) A 2-torus can be realized as a graph with 9 vertices. Its Betti vector is $(1,2,1)$.
   There is a Morse function which has 4 critical point, one with index $0$, two with index $1$ 
   and one with index $2$. \\
4) A Klein bottle has the Betti vector is $(1,1,0)$. Like the Torus, its Euler characteristic
   is zero but it is non-orientable ($b_2=0$). There is no Morse function with 2 critical points.
   But there is a Morse function with 4 critical points. \\
5) The Betti vector of a projective plane is $(1,0,0)$. There is no Morse function with 1 critical 
   point however. The minimal number of Morse critical points is 3. 
\index{2-sphere}
\index{3=sphere}
\index{Klein bottle}
\index{Projective plane.}
\index{Reeb theorem}

\paragraph{}
We first look what happens with the cohomology if $S$ is a sub-graph of a graph $G$ that is 
a $(k-1)$-sphere and we make an extension $G \to G +_S v$ along $S$. There
are two possibilities: either $S$ is the boundary of a contractible graph 
$B$ in $G$. In that case, we increase $b_k$ because we add a new $k$-sphere. 
If $S$ is not a boundary of a contractible part of $G$, 
it carries a cohomology (let the heat flow act to get a class). 
There are therefore two possibilities when attaching a $k$-ball along 
a $(k-1)$-sphere: 

\begin{lemma}
If a $k$ ball $B = S+x$ (handle) is attached to a $(k-1)$ sphere $S$, then either 
$b_k$ increases by $1$ or $b_{k-1}$ decreases by $1$.
\end{lemma}

\begin{thm}[Strong Morse inequalities]
For any Morse function $f$ on a graph, we have
$b_k-b_{k-1} + b_{k-2}- \cdots + (-1)^k b_0 \leq
 c_k-c_{k-1} + c_{k-2}- \cdots + (-1)^k c_0$.
\end{thm}
\begin{proof}
For each fixed $k$, we use induction with respect to the number $n$ of vertices in the graph. 
For $k=0$, the inequality $b_0 \leq c_0$ could also be proven directly. 
The left hand side is the number of connectivity components the right the 
minimal number of minima which a Morse function can have. \\
Lets now look at the induction step with respect to the number of vertices. 
When adding a new vertex, $b_0$ can only change if the number of connected components
changes. But then also $c_0$ changes. In general, $c_0$ the number of
local minima can increase in the Morse build-up. \\
For $k=1$, we have the inequality $b_0-b_1 \leq c_0-c_1$. It again is clear for
one vertex. If we look at the Morse build-up for $f$, the numbers can change
only at critical values. If a vertex $v$ with Morse index $k$ is added, the Euler
characteristic gets augmented by $(-1)^k$. Because we add a $(k+1)$-ball (handle)
along a $k$-sphere, the cohomology changes. \\
The Euler-Poincar\'e formula and Poincar\'e-Hopf show that $b_k$ can increase by $1$
if $c_k$ increases by $1$ or that $b_{k-1}$ decreases by $1$. 
The later happens if we cover an existing open $(k-1)$-sphere with a $k$ ball. 
The former happens if we cover an existing
$k$-ball along the boundary sphere with a $k$-ball, forming a new non-contractible
$k$-sphere. (see the lemma). In both cases, the inequalities remain valid. 
It is the second possibility which is the reason why we do not have equality in general.
Every time we remove a cohomology class, the total (unsigned) sum $c=\sum_k c_k$ increases,
while the total Betti number (unsighed) $\sum_k b_k$ decreases.
\end{proof}
\index{Morse inequality}
\index{Morse inequality ! strong}

\paragraph{}
By adding successive equations in this list, we get the {\bf weak Morse inequalities}
$b_k \leq c_k$. An other consequence is obtained for $k=d$, we have equality by 
the Poincar\'e-Hopf formula:  $b_d-b_{d-1} + b_{d-2}- \cdots + (-1)^d b_0 = 
 c_d-c_{d-1} + c_{d-2}- \cdots + (-1)^d c_0$ because both sides are then
the Euler characteristic. 
\index{weak Morse inequality}
\index{Morse inequality ! weak}

\paragraph{}
{\bf Examples:}  \\
1) For the 1-torus (circle), we have
$b_0 =1$ , $b_1=1$, $c_0=1$ , $c_1=1$.
Lets take $G=C_4 = \{ \{1\},\{2\},\{3\},\{4\},\{1,2\},\{1,3\},\{3,4\},\{4,2\} \}$
where we numbered the points so that $g(v)=v$ is a Morse function. 
Start with building up the graph $\{ \{1\} \}$ which is a critical point of 
Morse index $0$. When adding the points $2$ and $3$ we have homotopy extensions where
$S^-g(v)$ is $K_1$ and so contractible. When adding the 4th point, we add a $1$-handle
to a $0$-sphere $S=\{ \{2\},\{3\} \}$.  \\
2) If we take the $C_4$ from above and add an other 1-ball to $S=\{\{2\},\{3\}$,
then $b_1$ is increased by $1$ and the Betti vector becomes $(1,2)$. \\
3) We we take the $G=C_4$ from above an add a $2$-ball to the $1$-sphere $S=G$
then $b_1$ decreases by $1$.  \\
4) For the 2-torus we have  $b_0 =1$ , $b_1=2$ , 
$b_2=1$, $c_0=1$ , $c_1=2$ , $c_2=1$. The minimal number of critical points which 
a Morse function can have is $4$. 

\section{Category Theorem} 

\paragraph{}
In this section, we look at the fundamental theorem of Lusternik-Schnirelmann theory. 
\cite{CLOT} write: 
{\it LS category is like a Picasso painting. Looking at category from different perspectives
produces completely different impressions of category's beauty and applicability.}. 

\paragraph{}
Recall that ${\rm cri}(G)$ is the minimal number of critical points which a 
locally injective function $g$ on 
$G$ can have. This means that the function $g$ does not necessarily have to be a Morse function. 
And ${\rm cat}(G)$ is the minimal number of contractible graphs which cover 
the graph, which means both vertices and edges of course as a spanning tree already covers
all vertices. First to the lower bound: 

\begin{thm} 
${\rm cat}(G) \leq {\rm cri}(G)$.
\end{thm}

\begin{proof}
We use induction with respect to the minimal number $k$ of critical points of $G$.
Without loss of generality, we can assume that the function $g$ is injective. 
The reason is that a general locally injective function can be modified to be 
injective without changing the critical points. The function values of $g$
now produce a total order on the vertices $V$ and so a stratification, 
where $G_j$ is the sub-graph of $G$ generated by the vertices $\{v_1, \dots, v_j\}$. \\
Let us start with the induction foundation: 
if $k=1$, then by definition, $G$ is contractible, where the minimum is the only 
critical point. Because no further critical point occurs, all other extensions
$G_j \to G_{j+1}$ are homotopy steps keeping $G_{j+1}$ contractible if $G_j$ is
contractible. \\
Now lets look at the induction step. Lets assume that every graph with maximally 
$k$ critical points can be covered with $k$ contractible graphs. Now take a graph 
with $k+1$ critical points and that $v_{m_k}$ is the largest critical point.
Now all $G_{m_k+j}$ are contractible until we reach $G_{m_{k+1}-1}$. The next
expansion gives us $v_{m_{k+1}}$, the last critical point. 
As $G_{m_{k+1}-1}$ is contractible and the unit ball $B(v_{k+1})$ 
is contractible we could get over the $k+1$'th critical point by adding an other
contractible graph. Now, there are no critical points any more and the rest of the
extensions $G_j$ for $j=m_{k+1}+1$ to the last vertex $n$ are all homotopy extensions
not changing the number of contractible parts.
\end{proof} 

\paragraph{}
One could also argue closer to the continuum: 
define the set $S_k$ of sub-graphs of $G$ of category 
larger or equal than $k$ in M is not empty. We can assume that each graph in
$S_k$ is maximal in the sense that it is not a subgraph of a larger graph with 
the same category $k$. Take any enumeration function $g$.
The value $c_k = \min_{S \in S_k} \max_{x \in S} g(x)$ is a critical value
because if it were not, $S$ would not be maximal and we could extend $S$ without
changing category. It follows that the number of critical points is 
bigger or equal than ${\rm cat}(M)$. 

\paragraph{}
It follows that if there is a function with only one critical point, then it
produces a contractible graph. Lets move to the lower cohomological bound: 

\begin{thm}
${\rm cup}(G)+1 \leq  {\rm cat}(G)$.
\end{thm}

\paragraph{}
First a lemma: 

\begin{lemma}
Given any $p \geq 1$-form $f$ and any contractible 
sub-graph $U$ of $G$, there is a coboundary $h=dg$ such that $f-h$ is 
zero on $U$. 
\end{lemma}

\begin{proof} 
Because the cohomology of $U$ is trivial, every 
cocycle is a coboundary. That is, $f$ restricted to $U$ is a coboundary
when we look at it as a k-form on $U$. This means $f|U=dg$,where $g|U$ means $g$
restricted to $U$.  So, $f-h=f-dg$ is zero on $U$.
\end{proof}

\paragraph{}
Now to the proof of the lower bound through cup length: 

\begin{proof} 
Assume that ${\rm cup}(G)=k$. 
Take a maximal set of forms $f_1, \dots ,f_k$, where $f_j$ is a $p_j$-form
with $p_j \geq 1$ and such that $f_1 \wedge f_2 \wedge \cdots \wedge f_k \neq 0$. 
We want to show that ${\rm cap}(G) \geq k+1$. Assume this is not true, then 
we could find sets $U_1, \dots,U_k$ that are contractible and cover $G$. 
By the above lemma we can deform each $f_k$ by adding coboundaries 
such that $f_k \in H^1(G,U_k)$, meaning that that $f_k$ are zero in $U_k$.
But this means that the cup product of the $f_k$ is zero. 
\end{proof}

\paragraph{}
The following theorem is adapted from the continuum, where it applies to 
manifolds.  
Remember that {\bf order} of an open cover is the least integer $k$ such that 
there are $k+1$ members of the covering with non-zero intersection. 
The {\bf topological dimension} was the smallest $k \geq 0$ such that for
every open cover $\mathcal{U}$ of $G$, there is a refinement of order $k$. 
We have seen that the topological dimension of a graph 
agrees with the maximal dimension of the graph. 

\begin{thm}
If $G$ is a connected manifold then
${\rm cat}(G) \leq {\rm dim}(G)+1$.
\end{thm}
\begin{proof}
Let us denote the dimension with $d$. In order to prove the inequality
we have to construct $(d+1)$ contractible graphs which cover $G$.
We use induction with respect to $k$. Let us look at the case $k=1$.
It can be covered with two linear graphs.
As for the induction step:
we can cover each unit sphere $S(x)$ with $d$ contractible
graphs of dimension one less. Let $U(x,i)$ denote these graphs. 
Now build a spanning tree $T_i$ connecting all $\{ U(v,i),  v \in V \}$
and form the $d$ graphs $U_i = T_i \cup \{ U(v,i),  v \in V \}$. These graphs are 
contractible and their union cover all unit spheres simultaneously. 
Additionally take an other spanning tree $U_0$, hitting all vertices. Also this is 
contractible. The union $U_0 \cup \{ U_i \}_{i =1}^d$ covers the entire graph. 
\end{proof} 

\paragraph{}
The proof does not extend to general graphs. Already the induction foundation
fails. For a one dimensional graph, the category is the {\bf edge arboricity},
the minimal number of trees covering the graph. By the Nash-Williams theorem, this is
bounded below by the smallest integer larger or equal than $(|E|-1)/|V|$.
Lets take $K_n$ with $n$ vertices and $n(n-1)/2$ edges. Make a one dimensional 
refinement. Then we have $n+n(n-1)/2$ vertices and $n(n-1)$ edges.
For one dimensional graphs without triangles, the arboricity and so 
category can already become arbitrarily large. 
\index{arboricity}
\index{Nash-Williams theorem}

\paragraph{}
As mentioned, category is in graph theory close to{\bf arboricity}.
The {\bf Nash-Williams theorem} gives the arboricity in terms of the maximal 
$|E_H|/(|V_H|-1)$ taken over the set of all {\bf induced sub-graphs} $H$ of $G$. 
We can also look at coverings by edges which is an {\bf edge coloring problem}. 
The {\bf edge chromatic number} or {\bf edge index} is by Vizing's theorem 
between $\Delta(G)$ and $\Delta(G)+1$, where $\Delta(G)$ is the maximal vertex degree of the graph. 
We have worked on arboricity and chromatic number questions elsewhere 
\cite{ArboricityManifolds,ThreeTreeTheorem}.
\index{arboricity} 
\index{vertex degree}
\index{chromatic index} 
\index{Vizing theorem}

\paragraph{}
There are some higher dimensional discrete manifolds, for which we can compute the 
Lusternik-Schnirelmann category.  Since ${\rm cup}(G)=n$ for a $n$-dimensional torus, we have 
${\rm cat}(G) \geq n+1$. 
We need therefore at least $n+1$ contractible sub-graphs to cover $G$.
To get the upper bound construct a function with $n+1$ critical points.
We especially have ${\rm cat}(T^2) \geq 3$. 
On the 2-torus one can get $3$ critical points so that ${\rm cri}(G) \leq 3$. 
The book \cite{CLOT} gives the example $f(x,y)=\sin(x) \sin(y) \sin(x+y)$ which 
has only 3 critical points. Can we explicitly construct a function with $n+1$ 
critical points on a discrete $n$ torus? 

\paragraph{}
Every tree, a connected simply connected graph has category $1$. A connected 
graph $G=(V,E)$ without any triangles defines a one-dimensional 
simplicial complex $G = \{ \{v\},v \in V \} \cup \{ \{(a,b),(a,b) \in E \}$.
The Lusternik-Schnirelmann category of any such graph is $1$ if it is simply 
connected and larger if it is not simply connected. Since the maximal dimension of $G$
is $1$, the cup -length is $0$ or $1$. 
What is the minimal number of critical points? 
For a small figure 8 graph, we can just brute force
all functions and see that there are always at least 3 critical points. 

\paragraph{}
It would be nice to get a {\bf Morse cohomology} in the discrete. There
is hope: if $G$ is a graph and $G_1$ is the Barycentric refinement
in which the complete sub-graphs are the vertices and two are connected
if one is contained in the other, then $S^-(x)$ for
$f(x) = {\rm dim}(x)$ is a Morse function on $G_1$ because $S_f^-(x)$ is
a $k-1$-dimensional sphere if the dimension of $x$ was $k$. This means that 
the Morse index is $k$ and $i_g(x)=(-1)^{{\rm ind}(x)} = \omega(x)$. 
So, $\chi(G_1)=\sum_x \omega(x)$ by Poincar\'e Hopf. But this is the definition
of $\chi(G)$. We see $\chi(G)=\chi(G_1)$. 
In this case, $S^-(x) + S^+(x) = S(x)$, where $A+B$
is the {\bf join} of two graphs $A,B$ defined as 
$A+B=(V(A) \cup V(B),E(A) \cup E(B) \cup \{ (a,b),a \in A,b \in B \})$. 
Every $v \in G_1$ is a critical point. The Morse cohomology boundary agrees
now with the usual exterior derivative. 
$d(x,y)={\rm sign}(x,y)$ which is $1$ if $y \subset x$ has dimension $1$ less,
and matches orientation, $-1$ if $y \subset x$ has dimension $1$ less
and does not match orientation and $0$ else. 

\paragraph{}
Here is a nice example showing that counting is a Morse buildup 
\cite{Bjoerner2011,CountingAndCohomology} \footnote{We found \cite{Bjoerner2011} 
only after writing \cite{CountingAndCohomology}}.
The set of integers $\{2,3, \dots, n\}$ define a graph in which 
two nodes=integers are connected if one divides the other. We can see this graph as the 
Barycentric refinement of the {\bf complete prime graph} $P$, 
Adding a new number is
a Morse extension. Critical points are integers $n$ that have no square prime factors.
The Morse index is then the value of the Moebius function $-\mu(n)$. 
It follows from the Poincar\'e-Hopf formula that
$\chi(G) = 1-M(n)$ is the Euler characteristic, where $M(n)$ is the 
{\bf Mertens function} $M(n) = \sum_{k=1}^n \mu(k)$. 
\index{Mertens function}
\index{prime graph}
\index{Moebius function} 

\paragraph{}
We should point out {\bf some variations} about Lusternik-Schnirelmann theory in the
literature: classically, one sometimes uses the {\bf reduced category} ${\rm cat}(A)-1$ or
one then one uses uses ${\rm cup}(G) + 1$ and calls this {\bf augmented cup length} the 
{\bf algebraic category}.
Classically, Lusternik-Schnirelmann category is often silently defined only for connected spaces and
with respect to an ambient space in that one defines ${\rm cat}_G(A)$. 
The inequality ${\rm cat}(G) \leq {\rm dim}(G)+1$
for example only can hold for connected spaces. Example:
take a zero-dimensional space with $k$ connected components.
This space has category $k$ but is zero dimensional so that ${\rm cat}(G) \leq {\rm dim}(G)+1$
definitely can not hold. As for contractibility with respect
to an ambient space or contractibility independent of a background, one classically defines
${\rm cat}(A)=1$ if the injection $i:A \to X$ is
{\bf homotopic to a constant}. A circle $A$ in a sphere $G$ with that definition has 
category ${\rm cat}_G(A)=1$, while a circle ${\rm cat}_A(A)=2$ in itself has category $2$. 
The relative category ${\bf cat}_G(A)$ can also be defined for graphs by declaring it to be
the minimal number of contractible sets $U_i$ in $G$ which cover $A$. 
\index{Relative category}
\index{reduced category}
\index{augmented cup length}

\paragraph{}
In order to push the dimension result ${\rm cat}(A) \leq {\rm dim}(A) +1$ to the discrete, we
need a topology for which  connectedness match the definitions for graphs. 
In 2017 we defined a Zariski type topology on the simplex set $G$ where the {\bf closed sets}
are the simplicial complexes of sub-graphs $A$ of $G$. 
The open sets includes the {\bf void} $\emptyset$ which is also the simplicial complex 
of the empty graph and so also closed. 
\footnote{A simplicial complex does not contain 
the empty set. The axiom applies for the void $\emptyset$ as it does not contain the empty set. }
We here defined the graph topology as the topology with basis $B^+(x) = S^+(x) \cup \{x\}$ with 
$x \in G$. Closing it under intersections and union produces a {\bf topology} 
on the set of simplices. The closed sets are generated by the basis $B^-(x) = S^-(x) \cup \{x\}$.
In a cyclic graph $C_n$ for example, a set $S^+(\{b\}) = \{ \{b\},\{a,b\},\{b,c\} \}$ is
an example of an open set as is $\{ \{a,b\} \}$ the intersection of two
such $B^+(x) \cap B^+(y)$ with $(x,y) \in E$. The set 
$S^-(\{a,b\}) = \{ \{ a\}, \{b\} ,\{a,b\} \}$ is a closed set. The open basis 
elements play the role of open intervals the closed basis elements the role of closed intervals. 
Two different open sets intersect in an open set and two closed sets intersect in a closed set.
\index{void}
\index{Zariski type}

\paragraph{}
Let us take a locally injective function $f$ and a critical point $v$ for $f$ on 
a general finite simple graph $G$. Let us assume that the index $i_g(v)$ of the
vertex $v$ is zero. This implies that $S^-_f(v)$ has Euler characteristic $1$. 
Sometimes we can modify $f$ locally without changing the nature of other critical 
points so that $v$ becomes a regular point, sometimes we have an 
essential singularities and
the singularity is not a removable singularities.  Let us take the example of a large
graph $G$ which contains a vertex $v$ on which $S(v)$ is a {\bf dunce hat} and such that
within $G$, we can make a homotopy deformation of $S(v)$ to a contractible graph. 
If we assume that $f$ is a locally injective function on $G$ such that $v$ is a local
maximum then $v$ is by definition a critical point because $S(v)$ is not contractible.
But it is also a removable critical point because we can modify $f$ along the homotopy
deformation of $S(v)$ to a contractible space to a function $g$ without introducing new
critical points (we can assure this by making $G$ large enough) such that $g$ has now
$v$ as a regular point. In general this is not possible. Even if we take the cone $G=D+\{v\}$
over the dunce hat $D$ (meaning that $v$ is attached to $D$ such that $S(v)=D$ in $G$),
we have an essential singularity. An other example is the suspension $G$ of a discrete 
projective plane $P$. If $v$ is one of the two new vertices with $S(v)=P$ and $f$ is a
maximum on $v$ then $S_f^-(v)=S(v)=P$ but no other function on $G$ can render $S_f(v)$
contractible. 

\paragraph{}
The minimal length of a non-contractible closed simple path in a metric space is 
known as a {\bf systole}. \footnote{For 1-dimensional complexes, this is known as 
{\bf girth}. With that definition, if there is triangle present, the girth would be 3. 
For us, a triangle is contractible so that its systole is $0$, the systole 
of a d-sphere $0$ for $d>1$.} 
Gromov showed that for fixed $d$, all $d$-manifolds $M$ satisfy 
${\rm sys}^d(M) \leq C {\rm vol}(M)$ and showed also this is sharp as the 
Fubini-Study metric on $\mathbb{CP}^n$ gives $C=n!$. 
The {\bf systole category} is defined as the longest product of systoles which
give a curvature-free lower bound for the total curvature. This systole category
is a lower bound for Lusternik-Schnirelmann category. 
In the set-up considered here, the systole category is smaller or equal 
than the cup length. {\bf Gromov's systolic inequality} goes over to discrete $d$-manifolds.
One can take as {\bf volume} the number of maximal simplices. 
The discrete version of Gromov's inequality in graph theory would readily follow from 
the continuum by geometric realization. We are convinced that a finite 
proof would not be too difficult. For more information about the constant $C$ see
\cite{KatzSystolic}.
\index{systole}
\index{contractible path}
\index{fundamental group}
\index{girth}
\index{systole category}
\index{Fubini-Study metric}

\paragraph{}
Morse theory has grown in the last decades even more. 
{\bf Morse cohomology} is computationally often less complex than simplicial cohomology. 
The reason is that the differential complex is smaller and the matrices become smaller and more workable. 
Given a Morse function $f$, look at the vector space of all functions
on the set of critical points $v \in V$ of a Morse function. Define
$dg(p) = \sum_{q, {\rm ind}(q)={\rm ind}(p)-1} i_g(p,q) g(q)$, where for $p,q \in V$,
the index $i_f(p,q)$ is the signature of the intersection $W^+(p) \cap W^-(q)$.
It is defined as $\sum_{x \in W^+(p) \cap W^-(q)} {\rm sign}(x,W^+(p)) {\rm sign}(x,W^-(q)) \omega(x)$
and where ${\rm sign}(x,W^-(v))$ is $1$ if the orientation of $x$ agrees with the orientation of $q$
propagated along $W^-(q)$ and $(-1)$ else. One has to check $d^2=0$ to see
that there is a cohomology, called {\bf Morse cohomology}. 
In the case $g={\rm dim}$ on a Barycentric refinement, the Morse cohomology 
agrees with simplicial cohomology because $W^+(p) \cap W^-(q)= {\rm sign}(p,q)$.
\index{Morse cohomology}

\paragraph{}
The map which assigns to a signed graph $G$ in the Shannon ring 
its cohomology ring is a functor from the
category of signed graphs to the category of commutative rings.
Under addition in the Shannon ring, we get the direct product of rings.
Under multiplication in the Shannon ring, we get the tensor product of rings.
The functor therefore extends to a ring homomorphism from the ring of graphs
to a ring of cohomology rings of graphs.
We should add to \cite{ArithmeticGraphs, StrongRing,numbersandgraphs},
that the finite topos of delta sets can naturally be extended to an
associative, commutative ring with $1$. Introducing ``negative space"
can be done in canonical way. The disjoin union which is the
coproduct in the category defines a monoid and so extends naturally
to a group. Every ring element now can be written as $A-B$,
where $A,B$ are topos elements. The product
extends. Euler characteristic and cohomology extends, where the
Betti numbers just flip signs. Notions like
Poincar\'e-Hopf index or curvature extend.
\index{Shannon ring}
\index{arithmetic of graphs}

\printindex

\bibliographystyle{plain}

\end{document}